\documentclass[numbers,webpdf,imanum]{ima-authoring-template}%

\graphicspath{{img/}}

\theoremstyle{thmstyletwo}%
\newtheorem{theorem}{Theorem}
\newtheorem{proposition}[theorem]{Proposition}%

\newtheorem{definition}{Definition}

\usepackage{cancel}
\usepackage{relsize}
\usepackage{stmaryrd}
\usepackage{comment}
\usepackage{caption}
\usepackage{tikz}
\usepackage{booktabs}

\numberwithin{equation}{section}

\begin{document}

\DOI{DOI HERE}
\copyrightyear{2026}
\vol{00}
\pubyear{2026}
\access{Advance Access Publication Date: Day Month Year}
\appnotes{Paper}
\copyrightstatement{Published by Oxford University Press on behalf of the Institute of Mathematics and its Applications. All rights reserved.}
\firstpage{1}


\title[Global Iterative Methods for SPAIs of SPD matrices]{Global Iterative Methods for Sparse Approximate Inverses of Symmetric Positive Definite Matrices}

\author{Nicolas Venkovic and Hartwig Anzt
\address{\orgdiv{Chair of Computational Mathematics\\ School of Computation, Information and Technology (CIT)}\\ \orgname{Technical University of Munich}\\\orgaddress{Heilbronn, \postcode{74076}, \country{Germany}}}}

\authormark{Venkovic and Anzt}

\corresp[*]{Corresponding author: \href{email:nicolas.venkovic@tum.de}{nicolas.venkovic@tum.de,venkovic@gmail.com}}

\received{20}{08}{2026}
\revised{Date}{0}{Year}
\accepted{Date}{0}{Year}

\editor{Associate Editor: Professor Endre Suli}

\abstract{This work is motivated by symmetric positive definite (SPD) matrices for which the best sparse approximate inverse (SPAI) with the prescribed nonzero pattern of $A^k$ for some moderate value of $k$, e.g., 1, 2, 3, or 4, fails to capture essential features of the inverse, such as definiteness.
In this context, we consider short-recurrence iterative methods for the computation of SPAIs, that is, methods with sparse matrix iterates whose nonzero structure is globally updated at each iteration based on short-recurrence relations.
In particular, we consider the minimal residual (MR) method, its newly proposed variant enriched with one previous search direction, namely the locally optimal minimal residual (LOMR) method, and the conjugate gradient (CG) method with sparse matrix iterates and Frobenius inner products.
We show that, for SPD matrices, the MR method converges linearly with rate $(1-\lambda_{\min}^2/\operatorname{tr}(A^2))^{1/2}$, irrespective of the initial guess. 
While LOMR inherits unconditional monotone convergence from MR, its observed convergence behavior is that of a monotonically decreasing lower envelope to the CG residual norm without the occasional spurious oscillations proper to CG.
All three methods are implemented with practical dropping strategies to control the growth of nonzero patterns in the approximate inverse.
Numerical experiments are performed where SPAIs are computed with a prescribed cap on density, and the performance of those SPAIs as preconditioners for CG solves is quantitatively assessed for each method.
}

\maketitle

\section{Introduction}\label{sec:intro}
Sparse approximate inverses (SPAIs) have been used to precondition iterative solvers since the 1970s (see \cite{benson1973,frederickson1975fast}), in which case a sparse matrix $M\in\mathbb{R}^{n\times n}$ is formed explicitly to approximate the inverse $A^{-1}$ of some typically sparse matrix $A\in\mathbb{R}^{n\times n}$, and used as a preconditioner to accelerate the iterative solve of linear systems of the form $Ax=b$ (resp., $AX=B$), see~\cite{benzi1999comparative,benzi2002preconditioning} for overviews.
The main advantage of resorting to SPAIs is that their application, which boils down to the execution of sparse matrix-vector (resp., sparse matrix-matrix) kernels, is both portable and parallelizable.
The resort to SPAIs for preconditioning purposes relies on the assumption that sufficiently good but sparse approximations of $A^{-1}$ exist and can be found.
Strictly speaking, the inverse of a sparse matrix is generally dense.
That being said, it is rather common that many entries in the inverse of a sparse matrix have small magnitude, and dropping those values, that is, setting them to zero, minimally degrades the approximation of the inverse.
A particularly telling instance of this phenomenon is that of banded symmetric positive definite (SPD) matrices whose inverses have entries bounded by an exponentially decaying envelope along rows and columns~\cite{demko1984decay}.
Similarly, strongly diagonally dominant matrices also exhibit a rapid decay in the magnitude of the nonzero entries of their inverses.
This common feature of inverses of sparse matrices is leveraged for the design of SPAIs.

\begin{figure}[H]
\centering
\hspace*{-1.8cm}\includegraphics[scale=.8]{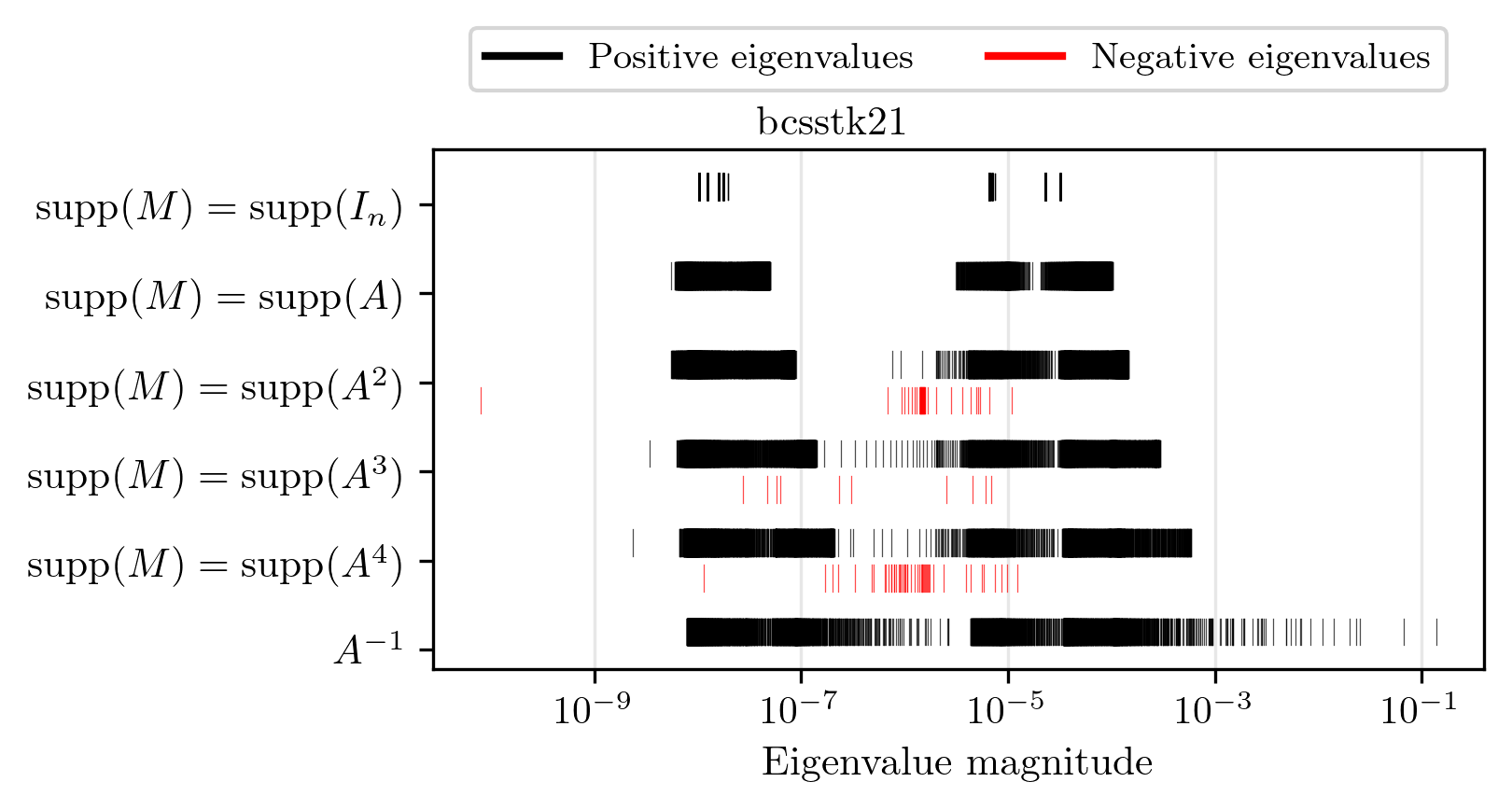}
\hspace*{1cm}\includegraphics[scale=.87]{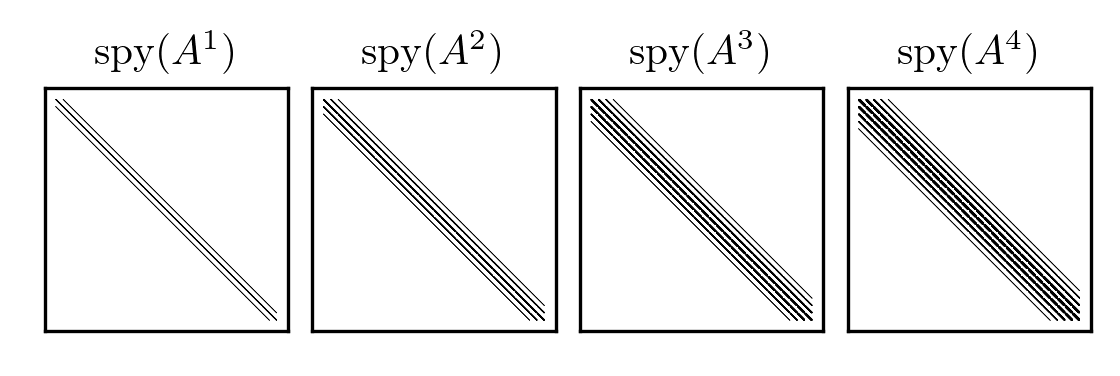}
\caption{Spectra of best SPAI with low-monomial sparsity for matrix bcsstk21 (Experiment09).}
\label{fig:lack-of-definiteness-1}
\end{figure}

\begin{figure}[H]
\centering
\hspace*{-1.8cm}\includegraphics[scale=.8]{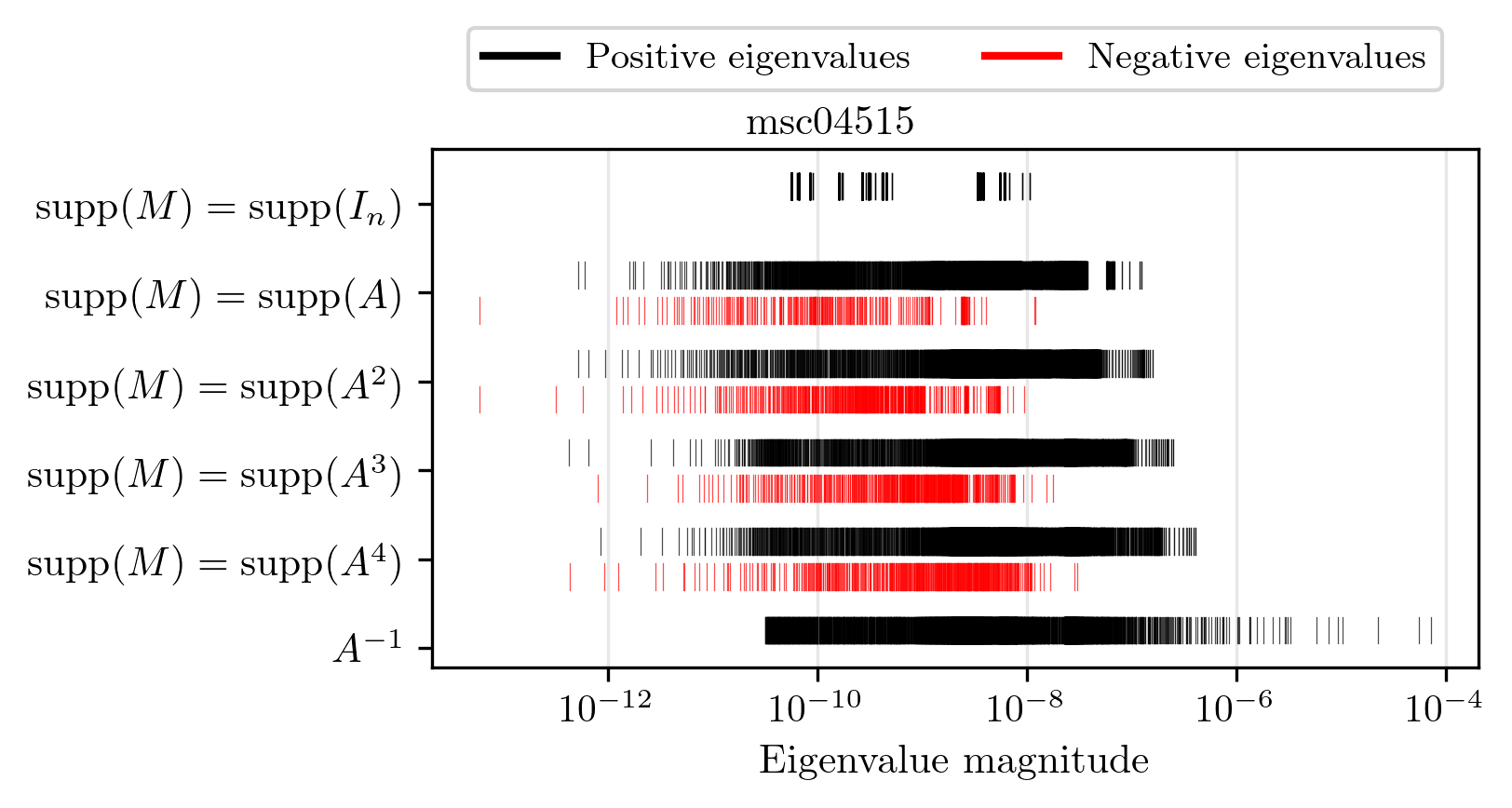}
\hspace*{1cm}\includegraphics[scale=.87]{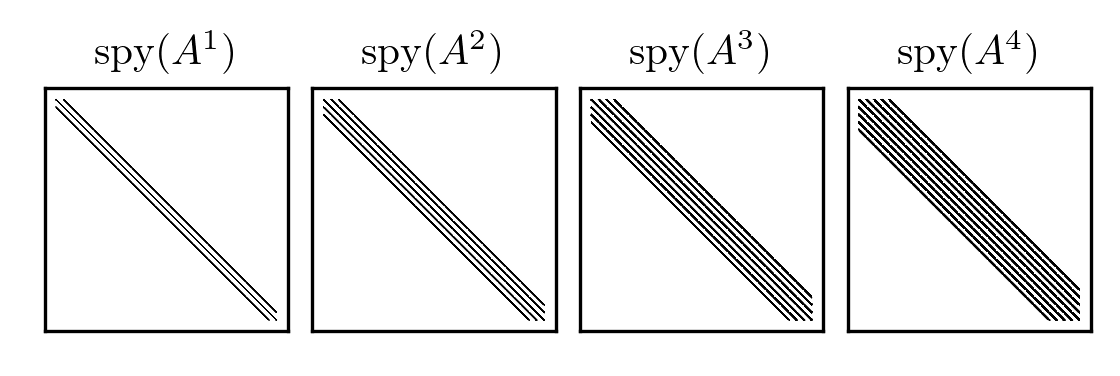}
\caption{Spectra of best SPAI with low-monomial sparsity for matrix msc04515 (Experiment09).}
\label{fig:lack-of-definiteness-2}
\end{figure}

Over time, several methods have been proposed for the computation of SPAIs, see~\cite{chow1997parallel,benzi1999comparative,Saad2003SparseBook,Bai2021MatrixAA} for overviews.
Perhaps the most significant distinction between those methods is related to whether the SPAI is expressed as a product of structured matrices, that is, a factorization, commonly referred to as a factorized sparse approximate inverse (FSAI), see~\cite{kolotilina1993factorized,kaporin1994new,benzi1996sparse,benzi1998sparse}, or simply as a sparse matrix, herein referred to as standard SPAI.
The advantage of FSAIs is that they enable more explicit control of spectral characteristics such as non-singularity and positive definiteness.
In this work, however, we focus on standard SPAIs and defer the treatment of FSAIs to another work.
Most existing approaches to compute SPAIs are iterative methods aimed at minimizing the residual Frobenius norm, one exception being~\cite{holland2005sparse}, where other norms are investigated.
At each iteration of an iterative method for the computation of SPAIs, an iterate $M_{i+1}$ of SPAI is formed from a previous iterate $M_i$ with the goal of minimizing $\|I_n-AM_{i+1}\|$ over a given approximation space.
The first documented attempts to develop such methods can be found in~\cite{benson1973,frederickson1975fast,benson1982iterative,benson1984parallel}.
In most of those early efforts, the sparsity pattern of the SPAI is assumed beforehand, eventually leading to decoupled and highly parallelizable least-squares problems.
For those approaches, common choices of nonzero patterns are extracted from small powers of $A$ which, in practice, can yield high nonzero densities, without guarantee of achieving a good approximation of $A^{-1}$.
Indeed, considering for example the two matrices bcsstk21 and msc04515 of the SuiteSparse Matrix Collection (\url{https://sparse.tamu.edu/}), see Figs.~\ref{fig:lack-of-definiteness-1} and \ref{fig:lack-of-definiteness-2} where the support denoted by $\operatorname{supp}(X)$ refers to the set of indices of nonzero elements of $X\in\mathbb{R}^{n\times n}$, we see that using nonzero patterns from $A$ to $A^4$ to construct SPAIs consistently fails to account for the positiveness of the inverse.
Such indefinite SPAIs of SPD matrices cannot be used as preconditioners of conjugate gradient (CG) iterations.

One of the most influential attempts to move away from the pre-definition of sparsity patterns is found in~\cite{Chow1998ApproximateIP}, where global iterative steepest descent (SD) and minimal residual (MR) methods are introduced, implemented with sparse data structures and computations, and deployed with dropping strategies, leading to nonzero patterns that evolve from one iteration to another with a contained nonzero density.
As of today, iterative methods based on the MR-based approaches of~\cite{Chow1998ApproximateIP} have become one of the most standard approaches to compute SPAIs.
As evidenced by experiments in Section~\ref{sec:dropping-free-experiments}, even when no dropping strategy is applied, global iterations based on the minimal residual and steepest descent methods often exhibit very slow convergence behaviors.
More particularly, when applied to SPD matrices, these methods consistently fail to produce SPD approximate inverses, rendering the SPAI unfit to precondition iterative solves by CG. 

In this work, we attempt to remedy these issues.
To do so, we introduce the locally optimal minimal residual (LOMR) method by enriching the search space with the previous search direction, and applying orthogonality constraints accordingly, as a means to accelerate convergence.
The LOMR method, that can be seen as a correction of MR iteration, is proven to converge globally, that is, independently of the initial iterate. 
In practice, we show experimental results that suggest that LOMR performs at least as well as sparse-mode CG iterations with matrix iterates all while showing monotone decrease of residual norm thereby not being prone to spurious oscillations that occasionally happen with CG.

This work is organized as follows.
First, the MR and SD methods are presented in Section~\ref{sec:one-dim-proj}.
We prove global convergence of the MR method for SPD systems and characterize the convergence behavior.
Then, the CG method is presented in Section~\ref{sec:cg}, and the LOMR method in Section~\ref{sec:lo-variants}.
Preconditioned iterations are described in Section~\ref{sec:precond-iter}, and dropping strategies in Section~\ref{sec:dropping}.
In Section~\ref{sec:numerical-experiments}, we present the results of numerical experiments which can be reproduced using the Julia package developed by the authors and made available at:
\begin{center}
\href{https://github.com/venkovic/SparseApproximateInverses.jl}{https://github.com/venkovic/SparseApproximateInverses.jl}.
\end{center}
Our conclusions are presented in Section~\ref{sec:conclusion}.

\section{One-dimensional projection methods}\label{sec:one-dim-proj}
Let $A\in\mathbb{R}^{n\times n}$ be an SPD matrix.
We are interested in finding a right-approximate inverse of $A$, that is, $M\in\mathbb{R}^{n\times n}$ such that $I_n-AM$ is small in some sense.
For that, let us consider the Frobenius inner product given by:
\begin{align}
(X,Y)_F\in\mathbb{R}^{n\times n}\times\mathbb{R}^{n\times n}\mapsto \text{tr}(X^TY)
\end{align}
with the induced norm given by:
\begin{align}
X\in\mathbb{R}^{n\times n}\mapsto(X,X)_F^{1/2}=:\|X\|_F
\end{align}
where $\|\cdot\|_F$ denotes the Frobenius norm.
Formally, we are interested in finding some approximation $M$ to the solution of
\begin{align}\label{eq:pb-linear-system}
AX=I_n
\end{align}
or, alternatively, we seek to minimize the objective function $f:\mathbb{R}^{n\times n}\rightarrow[0,\infty)$ given by
\begin{align}\label{eq:pb-objective-function}
f(M):=\|I_n-AM\|_F^2,
\end{align}
that is, the Frobenius norm of the residual of Eq.~\eqref{eq:pb-linear-system} for an approximation $M$ of $X=A^{-1}$.
Although Eq.~\eqref{eq:pb-objective-function} famously decomposes into $n$ independent $2$-norms of vector residuals (e.g., see Section~10.5 in~\cite{Saad2003SparseBook}) which allows for embarrassingly parallelizable strategies, in this work, we are exclusively concerned with global iterations.
That is, the matrix $M$ is considered as a whole, instead of a set of $n$ unrelated columns.

In this section, we present one-dimensional projection methods aimed at minimizing Eq.~\eqref{eq:pb-objective-function}.
One particular set of such methods is that of descent algorithms, which we define as follows.
Given an iterate $M_i\in\mathbb{R}^{n\times n}$  of right-approximate inverse of $A$, a descent method defines a new iterate $M_{i+1}\in\mathbb{R}^{n\times n}$ as follows for a given search direction $P_i\in\mathbb{R}^{n\times n}$:
\begin{align}\label{eq:descent-proj}
M_{i+1}\in M_i+\text{span}\{P_i\}
\;\text{ s.t. }\;
R_{i+1}:=I_n-AM_{i+1}\perp A\,\text{span}\{P_i\}
\;\text{ for }\;
i=0,1,\dots
\end{align}
where the orthogonality is stated in terms of the Frobenius inner product.
As long as $AP_i\neq 0_{n\times n}$, the application of the Petrov-Galerkin condition of Eq.~\eqref{eq:descent-proj} leads to the following update formula:
\begin{align}\label{eq:iter}
M_{i+1}=M_i+\alpha_i P_i
\;\text{ where }\;
\alpha_i=\frac{(R_i,AP_i)_F}{\|AP_i\|_F^2}
\;\text{ for }\;
i=0,1,\dots.
\end{align}
A standard result is that the iterate $M_{i+1}$ is defined by Eqs.~\eqref{eq:descent-proj}-\eqref{eq:iter} if and only if the following holds:
\begin{align}\label{eq:descent-opt}
M_{i+1}=\arg\min_{M\in M_i+\text{span}\{P_i\}}\|I_n-AM\|_F.
\end{align}
Thus, in a sense, $M_{i+1}$ is optimal over the affine search space $M_i+\text{span}\{P_i\}$.
Given an initial right-approximate inverse $M_0\in\mathbb{R}^{n\times n}$ of $A$, the update formula in Eq.~\eqref{eq:iter} is at the heart of global descent iterations.
Different variants of global descent methods are obtained depending on how the search direction $P_i$ is defined.
In what follows, we briefly review two existing variants, namely the minimal residual algorithm and the steepest descent method.

\subsection{Minimal residual method}
The simplest global descent method is obtained by letting the search direction $P_i$ be given by the residual $R_i$:
\begin{align}
P_i:=R_i=I_n-AM_i.
\end{align}
The resulting method, referred to as the minimal residual (MR) algorithm, is described in Algo.~\ref{algo:mr}.

\begin{algorithm}[ht]
\caption{MR($A$, $M_0$, $\texttt{ExplicitResidualUpdate}$)}
\label{algo:mr}
\begin{algorithmic}[1]
\State $R_0:=I_n-AM_0$
\For{$i=0,1,\dots$}
\State $\alpha_i:=(R_i,AR_i)_F/\|AR_i\|_F^2$
\State $M_{i+1}:=M_{i}+\alpha_i R_{i}$
\If{$\texttt{ExplicitResidualUpdate}=\texttt{True}$}
\State $R_{i+1}:=I_n-AM_{i+1}$
\Else
\State $R_{i+1}:=R_i-\alpha_i AR_i$
\EndIf
\EndFor
\end{algorithmic}
\end{algorithm}

For SPD matrices, the MR iteration converges linearly to $A^{-1}$ irrespective of the initial guess $M_0$.

\begin{theorem}[Convergence of MR iterates]\label{theo:mr-convergence}
Irrespective of the initial guess $M_0$, the MR iterates $M_0,M_1,\dots$ converge to the unique solution $M=A^{-1}$ such that $AM=I_n$.
They do so linearly at a rate $(1-\lambda_{\min}(A)^2/\operatorname{tr}(A^2))^{1/2}$.
\end{theorem}

\begin{proof}[Proof of Theo.~\ref{theo:mr-convergence}]
For every new MR iterate $M_{i+1}$, the residual norm is given by
\begin{align*}
\|R_{i+1}\|_F^2
=(R_{i+1},R_{i+1})_F
=(R_{i+1},R_i)_F-\alpha_i(R_{i+1},AR_i)_F\\
\end{align*}
where, by definition, we have the orthogonality condition $R_{i+1}\perp A\operatorname{span}\{R_i\}$ so that
\begin{align*}
\|R_{i+1}\|_F^2
=(R_{i+1},R_i)_F
\end{align*}
and
\begin{align}
\|R_{i+1}\|_F^2
=(R_{i+1},R_i)_F
=(R_i-\alpha_i AR_i,R_i)_F
=\|R_i\|_F^2-\alpha_i(AR_i,R_i)_F.
\end{align}
Using the fact that 
\begin{align}
\alpha_i=\frac{(AR_i,R_i)_F}{(AR_i,AR_i)_F},
\end{align}
we get
\begin{align}
\|R_{i+1}\|_F^2
=&\,\|R_i\|_F^2-\frac{(AR_i,R_i)_F}{(AR_i,AR_i)_F}(AR_i,R_i)_F\\
=&\,\|R_i\|_F^2\left(1-\frac{(AR_i,R_i)_F}{(AR_i,AR_i)_F}\frac{(AR_i,R_i)_F}{\|R_i\|_F^2}\right)\\
=&\,\|R_i\|_F^2\left(1-\frac{(AR_i,R_i)_F^2}{\|AR_i\|_F^2}\frac{\|R_i\|_F^2}{\|R_i\|_F^4}\right)
\end{align}
where, by consistency of the Frobenius norm, we have $\|AR_i\|_F\leq\|A\|_F\|R_i\|_F$ such that 
\begin{align}
\|R_{i+1}\|_F^2
\leq&\,\|R_i\|_F^2\left(1-\frac{(AR_i,R_i)_F^2}{\|A\|_F^2\|R_i\|_F^2}\frac{\|R_i\|_F^2}{\|R_i\|_F^4}\right)
\end{align}
eventually leading to
\begin{align}\label{eq:1000}
\|R_{i+1}\|_F^2
\leq&\,\|R_i\|_F^2\left(1-\frac{(AR_i,R_i)_F^2}{\|A\|_F^2(R_i,R_i)_F^2}\right).
\end{align}
Let us then reformulate the quotient as
\begin{align}
\frac{(AR_i,R_i)_F}{(R_i,R_i)_F}
=\frac{\operatorname{vec}(AR_i)^T\operatorname{vec}(R_i)}{\operatorname{vec}(R_i)^T\operatorname{vec}(R_i)}
\end{align}
where $\operatorname{vec}(AR_i)=(I_n\otimes A)\operatorname{vec}(R_i)$ so that
\begin{align}\label{eq:1010}
\frac{(AR_i,R_i)_F}{(R_i,R_i)_F}
=\frac{((I_n\otimes A)\operatorname{vec}(R_i))^T\operatorname{vec}(R_i)}{\operatorname{vec}(R_i)^T\operatorname{vec}(R_i)}
=\frac{\operatorname{vec}(R_i)^T(I_n\otimes A)\operatorname{vec}(R_i)}{\operatorname{vec}(R_i)^T\operatorname{vec}(R_i)}.
\end{align}
The spectrum of $I_n\otimes A$ is the same as that of $A$.
Thus, since $A$ is SPD, we have 
\begin{align}\label{eq:1020}
\frac{\operatorname{vec}(R_i)^T(I_n\otimes A)\operatorname{vec}(R_i)}{\operatorname{vec}(R_i)^T\operatorname{vec}(R_i)}\geq \lambda_{\min}(I_n\otimes A)=\lambda_{\min}(A)>0.
\end{align}
Combining Eqs.~\eqref{eq:1000}, \eqref{eq:1010} and \eqref{eq:1020} with the fact that $\|A\|_F^2=\operatorname{tr}(A^2)$, we obtain
\begin{align}
\|R_{i+1}\|_F^2\leq\|R_i\|^2\left(1-\frac{\lambda_{\min}(A)^2}{\operatorname{tr}(A^2)}\right)
\end{align}
and thus
\begin{align}
\|R_{i+1}\|_F\leq C\|R_i\|
\;\text{ with }\;
C:=\left(1-\frac{\lambda_{\min}(A)^2}{\operatorname{tr}(A^2)}\right)^{1/2}<1. 
\end{align}

\end{proof}

\subsection{Steepest descent method}
As shown in Section 10.5.2 of~\cite{Saad2003SparseBook}, and given the symmetry of $A$, the gradient of the objective function Eq.~\eqref{eq:pb-objective-function} is given by
\begin{align}\label{eq:gradient}
\nabla_{M}f(M)=-2A^TR=-2AR
\;\text{ where }\;
R:=I_n-AM.
\end{align}
We denote the gradient direction by $G_i:=-AR_i$, and the steepest descent (SD) method consists of letting the search direction be opposite to the gradient direction:
\begin{align}
P_i:=-G_i
\;\text{ for }\;i=0,1,\dots.
\end{align}
The resulting procedure is given in Algo.~\ref{algo:steepest-descent}.

\begin{algorithm}[ht]
\caption{SD($A$, $M_0$, $\texttt{ExplicitResidualUpdate}$)}
\label{algo:steepest-descent}
\begin{algorithmic}[1]
\State $R_0:=I_n-AM_0$
\State $P_0:=AR_0$
\For{$i=0,1,\dots$}
\State $\alpha_0:=(R_i,AP_i)_F/\|AP_i\|_F^2$
\State $M_{i+1}:=M_{i}+\alpha_i P_{i}$
\If{$\texttt{ExplicitResidualUpdate}=\texttt{True}$}
\State $R_{i+1}:=I_n-AM_{i+1}$
\Else
\State $R_{i+1}:=R_i-\alpha_i AP_i$
\EndIf
\State $P_{i+1}:=AR_{i+1}$
\EndFor
\end{algorithmic}
\end{algorithm}

\section{Conjugate gradient method}\label{sec:cg}
In this Section, we present and summarize properties of the conjugate gradient (CG) method~\cite{hestenes_methods_1952} for the computation of SPAIs.
The CG method is defined as a matrix iteration procedure in Definition~\ref{def:cg}, with the corresponding iterates given in Proposition~\ref{prop:cg-iterates}, and the associated procedure described in Algo.~\ref{algo:cg}.

\begin{definition}[CG method]\label{def:cg}
Given an SPD matrix $A\in\mathbb{R}^{n\times n}$ with a right-approximate inverse $M_0\in\mathbb{R}^{n\times n}$, a sequence of CG iterates of right-approximate inverses of $A$ is defined by
\begin{align}\label{eq:700}
M_{i+1}\in M_{i}+\text{span}\{P_i\}
\;\text{ s.t. }\;
R_{i+1}:=I_n-AM_{i+1}\perp\text{span}\{P_i\}
\;\text{ for }\;
i=0,1,\dots,
\end{align}
where $P_i\in\mathbb{R}^{n\times n}$ is a search direction iterate defined as
\begin{align}\label{eq:705}
P_{i}\in R_i+\text{span}\{P_{i-1}\}
\;\text{ s.t. }\;
P_{i}\perp A\,\text{span}\{P_{i-1}\}
\;\text{ for }\;
i=1,2,\dots
\end{align}
with $P_0:=R_0$.
\end{definition}

\begin{proposition}[CG iterates]\label{prop:cg-iterates}
The iterates of the CG method (Definition~\ref{def:cg}) are given by
\begin{align}\label{eq:710}
M_{i+1}:=M_{i}+\alpha_iP_i
\;\text{ where }\;
\alpha_i
:=\frac{(R_i,R_i)_F}{(P_i,AP_i)_F}
\;\text{ for }\;
i=0,1,\dots,
\end{align}
in which the search direction is updated as follows:
\begin{align}\label{eq:715}
P_{i}:=R_i+\beta_iP_{i-1}
\;\text{ where }\;
\beta_i:=\frac{(R_{i},R_{i})_F}{(R_{i-1},R_{i-1})_F}
\;\text{ for }\;
i=1,2,\dots.
\end{align}
Then, the residuals are orthogonal.
That is:
\begin{align}\label{eq:720}
(R_{i+1},R_j)_F=0
\;\text{ for }\;
j=0,1,\dots,i.
\end{align}
The residuals are also orthogonal to the search directions:
\begin{align}\label{eq:725}
(R_{i+1},P_j)_F=0
\;\text{ for }\;
j=0,1,\dots,i,
\end{align}
and the search directions are $A$-orthogonal:
\begin{align}\label{eq:730}
(P_{i+1},AP_j)_F=0
\;\text{ for }\;
j=0,1,\dots,i.
\end{align}
\end{proposition}

\begin{proof}[Proof of Proposition~\ref{prop:cg-iterates}]
Beside the replacement of vector iterates by matrix iterates, and vector dot products by Frobenius inner products, the proof is state-of-the-art.
See, e.g., \cite{ciaramella_iterative_2022}.
\end{proof}

\begin{algorithm}[ht]
\caption{CG($A$, $M_0$, $\texttt{ExplicitResidualUpdate}$)}
\label{algo:cg}
\begin{algorithmic}[1]
\State $R_0:=I_n-AM_0$
\State $P_0=R_0$
\For{$i=0,1,\dots$}
\State $\alpha_i:=\|R_i\|_F^2/(P_i,AP_i)_F$
\State $M_{i+1}:=M_i+\alpha_i P_i$
\If{$\texttt{ExplicitResidualUpdate}=\texttt{True}$}
\State $R_{i+1}:=I_n-AM_{i+1}$
\Else
\State $R_{i+1}:=R_i-\alpha_iAP_i$
\EndIf
\State $\beta_{i+1}:=\|R_{i+1}\|_F^2/\|R_i\|_F^2$
\State $P_{i+1}:=R_{i+1}+\beta_{i+1} P_i$
\EndFor
\end{algorithmic}
\end{algorithm}

Although CG is formulated in Definition~\ref{def:cg} as a one-dimensional projection, it is globally optimal over a Krylov subspace of $A$ (see Proposition~\ref{prop:cg-krylov}) over which it minimizes the Frobenius $A$-norm of the error, as stated in Theorem~\ref{theo:cg-optimal}.
Additionally, the Frobenius $A$-norm of the error is bounded above as expressed in Eq.~\eqref{eq:732}.

\begin{proposition}[CG as a Krylov subspace method]\label{prop:cg-krylov}
The iterates of the CG method (Definition~\ref{def:cg}) are equivalently defined as
\begin{align}
M_i\in M_0+\mathcal{K}_i(A,R_0)
\;\text{ s.t. }\;
R_i:=I_n-AM_i\perp \mathcal{K}_i(A,R_0)
\;\text{ for }\;
i=1,2,\dots
\end{align}
where $\mathcal{K}_i(A,R_0)$ denotes the Krylov subspace of $A$ generated by the initial residual:
\begin{align}
\mathcal{K}_i(A,R_0):=\text{span}\{R_0,AR_0,\dots,A^{i-1}R_0\}.
\end{align}
Moreover, the search directions $P_0,\dots,P_{i-1}$ constitute an $A$-orthogonal basis of this Krylov subspace, which is equally spanned by the residuals:
\begin{align}
\text{span}\{P_0,P_1,\dots,P_{i-1}\}
=
\text{span}\{R_0,R_1,\dots,R_{i-1}\}
=
\mathcal{K}_i(A,R_0).
\end{align}
\end{proposition}

\begin{proof}[Proof of Prop.~\ref{prop:cg-krylov}]
The equivalence between one-dimensional iterates defined locally in Def.~\ref{def:cg} and the globally orthogonal Krylov subspace-based definition of CG is also covered in~\cite{ciaramella_iterative_2022} for vector iterates.
\end{proof}

\begin{theorem}[Optimality of CG iterates]\label{theo:cg-optimal}
The CG iterates are given by
\begin{align}
M_i\in M_0+\mathcal{K}_i(A,R_0)
\;\text{ s.t. }\;
R_i\perp\mathcal{K}(A,R_0)
\;\text{ for }\;
i=1,2,\dots
\end{align}
if and only if the right-approximate inverse $M_i$ is optimal in the sense that
\begin{align}
\|A^{-1}-M_i\|_{F,A}=
\underset{M_i\in M_0+\mathcal{K}_i(A,R_0)}{\min}\|A^{-1}-M_i\|_{F,A}
\;\text{ for }\;
i=1,2,\dots.
\end{align}
That is, the CG iterate minimizes the Frobenius $A$-norm of the error $E:=A^{-1}-M$ over the affine Krylov subspace of $A$ generated by the initial residual.
\end{theorem}
\begin{proof}[Proof of Theo.~\ref{theo:cg-optimal}]
See~\cite{Saad2003SparseBook} for case of vector iterates.
\end{proof}

\begin{theorem}[Convergence of CG iterates]\label{theo:cg-bound}
Let $A$ have $k\leq n$ distinct eigenvalues.
Then, the CG iterates (Definition~\ref{def:cg}) converge to $A^{-1}$ in no more than $k$ iterations.
Moreover, the Frobenius $A$-norm of the error is bounded as follows:
\begin{align}\label{eq:732}
\|A^{-1}-M_i\|_{F,A}
\leq
\sqrt{n}
\left(
\frac{\kappa(A)^{1/2}-1}{\kappa(A)^{1/2}+1}
\right)^i
\|A^{-1}-M_0\|_{F,A}
\;\text{ for }\;
i=0,1,\dots.
\end{align}
\end{theorem}
\begin{proof}[Proof of Theo.~\ref{theo:cg-bound}]
See~\cite{ciaramella_iterative_2022} for the case with vector iterates. Here, the factor $\sqrt{n}$ comes from the bounding of the Frobenius norm of the error polynomial of CG by the product of $\sqrt{n}$ with the 2-norm of the error matrix polynomial.
\end{proof}

\section{Locally optimal variants}\label{sec:lo-variants}
The concept of enriching the search space of optimal iterations with prior iteration history was introduced by Polyak~\cite{polyak_methods_1964} with a similar concept by Knyazev~\cite{knyazev2001toward} for the case of approximation of eigenpairs produced by CG \cite{knyazev1991preconditioned}.
The term "locally optimal" was coined by Knyazev~\cite{knyazev2001toward}.
Here, we consider so called locally optimal variants of the one-dimensional methods presented in this paper.
That is, the case where we minimize the objective function in Eq.~\eqref{eq:pb-objective-function}.
Standard descent algorithms such as MR and SD form a new iterate $M_{i+1}$ which achieves the minimum given by:
\begin{align}\label{eq:600}
\min_{M\in M_i+\text{span}\{P_i\}}\|I_n-AM\|_F.
\end{align}
We can think of the MR method as a trivial case of
\begin{align}\label{eq:602}
P_i\in\text{span}\{R_i,P_{i-1}\}.
\end{align}
Indeed, MR assumes $P_i=R_i$ for $i=0,1,\dots$.
Now, upon letting $P_{-1}:=0_{n\times n}$ and $P_0:=R_0$, we have that 
\begin{align}\label{eq:605}
\min_{M\in M_i+\text{span}\{R_i,P_{i-1}\}}\|I_n-AM\|_F
\leq
\min_{M\in M_i+\text{span}\{P_i\}}\|I_n-AM\|_F
<\|R_i\|_F
\;\text{ for }\;
i=0,1,\dots.
\end{align}
This observation that considering both $P_{i-1}$ and $R_i$ as simultaneous degrees of freedom to form a better search direction capable of achieving a smaller local minimum of residual norm $\|R_{i+1}\|_F$ is at the heart of the locally optimal minimal residual (LOMR) method we introduce in Definition~\ref{def:lomr}. 
The associated update formulas are given in Theorem~\ref{theo:lomr-iterates}, and the resulting algorithm in Algo.~\ref{algo:lomr}.
We can say that LOMR converges at least as fast MR, but in practice, experiments show that LOMR can converge significantly faster than MR.

\begin{definition}[LOMR method]\label{def:lomr}
Given an SPD matrix $A\in\mathbb{R}^{n\times n}$ with a right-approximate inverse $M_0\in\mathbb{R}^{n\times n}$, a sequence of LOMR iterates of right-approximate inverses of $A$ is defined by
\begin{align}\label{eq:610}
M_{i+1}:=
\arg\underset{M\in M_i+\text{span}\{R_i,P_{i-1}\}}{\min}\;
\|I_n-AM\|_F
\;\text{ for }\;
i=0,1,\dots,
\end{align}
where $P_{-1}:=0_{n\times n}$, $P_0:=R_0$ and $R_i:=I_n-AM_i$ for $i=0,1,\dots$.
\end{definition}

\begin{theorem}[]\label{theo:lomr-iterates}
The iterates of the LOMR method (Definition~\ref{def:lomr}) are given by
\begin{align}\label{eq:615}
M_{i+1}:=M_i+\delta_iR_i+\gamma_iP_{i-1}
\;\text{ for }\;
i=0,1,\dots
\end{align}
where
\begin{align}\label{eq:620}
\delta_i:=
\begin{cases}
(R_i,AR_i)_F/(AR_i,AR_i)_F&\text{if }i=0\\
[(AP_{i-1},AP_{i-1})_F\cdot(R_i,AR_i)_F-(AR_i,AP_{i-1})_F\cdot(R_i,AP_{i-1})_F]/c_i&\text{otherwise}\\
\end{cases}
\end{align}
and
\begin{align}\label{eq:625}
\gamma_i:=
\begin{cases}
1&\text{if }i=0\\
[(AR_i,AR_i)_F\cdot(R_i,AP_{i-1})_F-(AR_i,AP_{i-1})_F\cdot(R_i,AR_i)_F]/c_i&\text{otherwise}
\end{cases}
\end{align}
in which
\begin{align}\label{eq:630}
c_i:=(AR_i,AR_i)_F\cdot(AP_{i-1},AP_{i-1})_F-(AR_i,AP_{i-1})_F^2
\;\text{ for }\;i=1,2,\dots.
\end{align}
The residual $R_{i+1}$ and search direction $P_i$ are updated as follows:
\begin{align}\label{eq:635}
R_{i+1}:=&\,R_i-\delta_iAR_i-\gamma_iAP_{i-1},\\
P_i:=&\,\delta_iR_i+\gamma_iP_{i-1}.\label{eq:637}
\end{align}
\end{theorem}

\begin{proof}[Proof of Theorem~\ref{theo:lomr-iterates}]
The optimality condition of Eq.~\eqref{eq:605} may be recast as follows:
\begin{align*}
\|R_{i+1}\|_F=&\,
\underset{M\in M_i+\text{span}\{R_i,P_{i-1}\}}{\min}\;\|I_n-AM\|_F\\
=&\,
\underset{\Delta M\in\text{span}\{R_i,P_{i-1}\}}{\min}\;\|R_i-A\Delta M\|_F\\
=&\,\underset{\Delta R\in A\,\text{span}\{R_i,P_{i-1}\}}{\min}\;\|R_i-\Delta R\|_F.
\end{align*}
By the theorem of orthogonal projections, this last orthogonality statement holds if and only if
\begin{align}\label{eq:640}
\Delta R\in A\,\text{span}\{R_i,P_{i-1}\}
\;\text{ and }\;
R_i-\Delta R\perp A\,\text{span}\{R_i,P_{i-1}\},
\end{align}
so that there exist $\delta_i,\gamma_i\in\mathbb{R}$ such that
\begin{align}
\begin{cases}
(AR_i,R_i-A(\delta_iR_i+\gamma_i P_{i-1}))_F=0\\
(AP_{i-1},R_i-A(\delta_iR_i+\gamma_i P_{i-1}))_F=0\\
\end{cases}
\end{align}
which yields
\begin{align}
\begin{cases}
(AR_i,R_i)_F=\delta_i(AR_i,AR_i)_F+\gamma_i(AR_i,AP_{i-1})_F\\
(AP_{i-1},R_i)_F=\delta_i(AP_{i-1},AR_i)_F+\gamma_i(AP_{i-1},AP_{i-1})_F
\end{cases}
\end{align}
and whose solution is given by Eqs.~\eqref{eq:620}, \eqref{eq:625} and \eqref{eq:630} to form the 
update formula given by:
\begin{align*}
M_{i+1}:=M_i+\delta_iR_i+\gamma_iP_{i-1}
\;\text{ for }\;
i=1,2,\dots.
\end{align*}
For $i=0$, the update formula is the same as that of the MR method, that is:
\begin{align}
M_1:=M_0+\frac{(R_0,AR_0)_F}{(AR_0,AR_0)_F}R_0.
\end{align}
As for the search direction, the combination of Eqs.~\eqref{eq:602} and \eqref{eq:615} allows for different, but equally valid update formulae of the search direction.
In particular, we may assume
\begin{align}
M_{i+1}=M_i+P_i,
\end{align}
which implies the update formula given by Eq.~\eqref{eq:637}.
\end{proof}

\begin{algorithm}[ht]
\caption{LOMR($A$, $M_0$, $\texttt{ExplicitResidualUpdate}$)}
\label{algo:lomr}
\begin{algorithmic}[1]
\State $R_0:=I_n-AM_0$
\State $P_{-1}:=0_{n\times n}$
\For{$i=0,1,\dots$}
\If{$i=0$}
\State $\delta_i:=(R_i,AR_i,)_F/\|AR_i\|_F^2$
\State $\gamma_i:=1$
\Else
\State $c_i:=\|AR_i\|_F^2\cdot\|AP_{i-1}\|_F^2-(AR_i,AP_{i-1})_F^2$
\State $\delta_i:=\left.\left[\|AP_{i-1}\|_F^2\cdot(R_i,AR_i,)_F-(AR_i,AP_{i-1})_F\cdot(R_i,AP_{i-1})_F\right]\right/c_i$
\State $\gamma_i:=\left.\left[\|AR_i\|_F^2\cdot(R_i,AP_{i-1},)_F-(AR_i,AP_{i-1})_F\cdot(R_i,AR_i)_F\right]\right/c_i$
\EndIf
\State $M_{i+1}:=M_i+\delta_iR_i+\gamma_i P_{i-1}$
\If{$\texttt{ExplicitResidualUpdate}=\texttt{True}$}
\State $R_{i+1}:=I_n-AM_{i+1}$
\Else
\State $R_{i+1}:=R_i-\delta_iAR_i-\gamma_iAP_{i-1}$
\EndIf
\State $P_{i}:=\delta_iR_{i}+\gamma_iP_{i-1}$
\EndFor
\end{algorithmic}
\end{algorithm}

Note that, similarly to how we defined an LOMR method by enriching MR iterations, we could derive a locally optimal version of the SD method.
However, as we shall confirm in Section~\ref{sec:dropping-free-experiments}, the SD method offers poor convergence properties compared to MR, so that it is not worth defining an LOSD method.
As for the CG method, its global optimality over the affine span of search directions renders this method locally optimal in the sense stated in Eq.~\eqref{eq:605}.
Thus, if one were to define locally optimal versions of the CG method, they would essentially derive expensive procedures equivalent to the standard CG method.
If you wish to convince yourself of this, we invite the reader to experiment with the lopcg routine made available in the GitHub repository of this paper.

\section{Preconditioned iterations}\label{sec:precond-iter}
The convergence behavior of the iterative methods presented thus far is known to depend on the spectrum of $A$.
In practice, the number of iterations needed to achieve convergence may be so inconveniently large, that the use of a preconditioner becomes necessary.
For that matter, the simplest approach is perhaps that of left-preconditioning, that is, say we are equipped with a non-singular matrix $\Pi\in\mathbb{R}^{n\times n}$ which somehow approximates the inverse of $A$, and such that $X\in\mathbb{R}^{n\times n}\mapsto \Pi X$ can be computed efficiently.
Note that for the approach to remain consistent with our objective to construct a sparse approximate inverse, the preconditioner $\Pi$ needs also be sparse.
Then, instead of solving Eq.~\eqref{eq:pb-linear-system}, we search for approximate solutions of
\begin{align}\label{eq:500}
\Pi AX=\Pi
\end{align}
where the coefficient matrix $\Pi A$ is tentatively more amenable to efficient iterative solves.
Importantly, we note that, although $A$ is SPD, even if $\Pi$ is SPD, or simply symmetric, the coefficient matrix of Eq.~\eqref{eq:500} is generally not symmetric.
This lack of symmetry or positive-definiteness has to be accounted for for a proper application of preconditioning to matrix iterations.


Since all the iterations are defined here for a matrix $A$ self-adjoint with the inner product, this is leveraged as follows, as done commonly with CG vector iteration.
All the methods are defined with preconditioning applied from the left with the standard Frobenius inner product replaced by an inner product weighted by $\Pi^{-1}$, thus rendering $\Pi A$ self-adjoint with respect to the inner product.
Upon minimal derivative work, this eventually leads to the procedures given in Algos.~\ref{algo:pmr}-\ref{algo:lopmr}.

\begin{algorithm}[ht]
\caption{MR\_SPD($A$, $\Pi$, $M_0$, $\texttt{ExplicitResidualUpdate}$)}
\label{algo:pmr}
\begin{algorithmic}[1]
\Statex \Comment{$\Pi$ is an SPD preconditioner, e.g., $\Pi:=D_A^{-1}$}
\State $R_0:=I_n-AM_0$
\For{$i=0,1,\dots$}
\State $Z_i:=\Pi R_i$ 
\State $\alpha_i:=(Z_i,\Pi AZ_i)_F/(AZ_i,\Pi AZ)_F$
\State $M_{i+1}:=M_{i}+\alpha_i Z_{i}$
\If{$\texttt{ExplicitResidualUpdate}=\texttt{True}$}
\State $R_{i+1}:=I_n-AM_{i+1}$
\Else
\State $R_{i+1}:=R_i-\alpha_iA_iZ_i$
\EndIf
\EndFor
\end{algorithmic}
\end{algorithm}

\begin{algorithm}[ht]
\caption{SD\_SPD($A$, $\Pi$, $M_0$, $\texttt{ExplicitResidualUpdate}$)}
\label{algo:preconditioned-steepest-descent}
\begin{algorithmic}[1]
\Statex \Comment{$\Pi$ is an SPD preconditioner, e.g., $\Pi:=D_A^{-1}$}
\State $R_0:=I_n-AM_0$
\For{$i=0,1,\dots$}
\State{$Z_i:=\Pi R_i$}
\State $P_{i}:=AR_{i}$
\State $\alpha_i:=(Z_i,AP_i)_F/(AP_i,\|\Pi AP_i)F$
\State $M_{i+1}^{-1}:=M_{i}^{-1}+\alpha_i P_{i}$
\If{$\texttt{ExplicitResidualUpdate}=\texttt{True}$}
\State $R_{i+1}:=I_n-AM_{i+1}$
\Else
\State $R_{i+1}:=R_i-\alpha_iAP_i$
\EndIf
\EndFor
\end{algorithmic}
\end{algorithm}

\begin{algorithm}[ht]
\caption{CG($A$, $\Pi$, $M_0$, $\texttt{ExplicitResidualUpdate}$)}
\label{algo:pcg}
\begin{algorithmic}[1]
\Statex \Comment{$\Pi$ is an SPD preconditioner, e.g., $\Pi:=D_A^{-1}$}
\State $R_0:=I_n-AM_0$
\State $Z_0:=\Pi R_0$
\State $P_0=Z_0$
\For{$i=0,1,\dots$}
\State $\alpha_i:=(R_i,Z_i)_F/(P_i,AP_i)_F$
\State $M_{i+1}:=M_i+\alpha_i P_i$
\If{$\texttt{ExplicitResidualUpdate}=\texttt{True}$}
\State $R_{i+1}:=I_n-AM_{i+1}$
\Else
\State $R_{i+1}:=R_i-\alpha_iAP_i$
\EndIf
\State $Z_{i+1}:=\Pi R_{i+1}$
\State $\beta_{i+1}:=(R_{i+1},Z_{i+1})_F/(R_i,Z_i)_F$
\State $P_{i+1}:=Z_{i+1}+\beta_{i+1} P_i$
\EndFor
\end{algorithmic}
\end{algorithm}

\begin{algorithm}[ht]
\caption{LOMR\_SPD($A$, $\Pi$, $M_0$, $\texttt{ExplicitResidualUpdate}$)}
\label{algo:lopmr}
\begin{algorithmic}[1]
\Statex \Comment{$\Pi$ is an SPD preconditioner, e.g., $\Pi:=D_A^{-1}$}
\State $R_0:=I_n-AM_0$
\State $P_{-1}:=0_{n\times n}$
\State $Z_0:=\Pi R_0$
\For{$i=0,1,\dots$}
\If{$i=0$}
\State $\delta_i:=(Z_i,AZ_i,)_F/(AZ_i,\Pi AZ_i)_F$
\State $\gamma_i:=1$
\Else
\State $c_i:=(AZ_i,\Pi AZ_i)_F\cdot(AP_{i-1},\Pi AP_{i-1})_F-(AZ_i,\Pi AP_{i-1})_F^2$
\State $\delta_i:=\left.\left[(AP_{i-1},\Pi AP_{i-1})_F\cdot(Z_i,AZ_i,)_F-(AZ_i,\Pi AP_{i-1})_F\cdot(Z_i,AP_{i-1})_F\right]\right/c_i$
\State $\gamma_i:=\left.\left[(AZ_i,\Pi AZ_i)_F\cdot(Z_i,AP_{i-1},)_F-(AZ_i,\Pi AP_{i-1})_F\cdot(Z_i,AZ_i)_F\right]\right/c_i$
\EndIf
\State $M_{i+1}:=M_i+\delta_iZ_i+\gamma_i P_{i-1}$
\If{$\texttt{ExplicitResidualUpdate}=\texttt{True}$}
\State $R_{i+1}:=I_n-AM_{i+1}$
\Else
\State $R_{i+1}:=R_i-\delta_iAZ_i-\gamma_iAP_{i-1}$
\EndIf
\State $P_{i}:=\delta_i Z_{i}+\gamma_i P_{i-1}$
\State $Z_{i+1}:=\Pi R_{i+1}$
\EndFor
\end{algorithmic}
\end{algorithm}

\section{Dropping strategies}\label{sec:dropping}
As the iterative methods of Sections~\ref{sec:one-dim-proj} to \ref{sec:precond-iter} are deployed, the number of nonzero components in the main iterate $M$ and the search direction $P$, where we drop the iteration index for clarity, can grow uncontrollably, slowing down iterations, likely exceeding memory capacity.
To circumvent this issue, dropping strategies need be deployed.
That is, methods to set nonzero values in $M$ and $P$ to zero, thus relieving some of the burden of excessive storage requirements, all the while minimizing convergence hindering.
We present here hard thresholding and we refer the reader to Section~\ref{sec:heuristic} for a dropping heuristic where we draw inspiration from the work of~\cite{Chow1998ApproximateIP} with facts from when up to one nonzero element is dropped per column. 
However, we stress that, from our experience, sophisticated dropping strategy do not help the iteration more than hard thresholding.

\subsection{Hard thresholding}
The simpler method to drop nonzero components is called hard thresholding and may be applied both to the main iterate $M$ and the search direction $P$.
For example, let the nonzero components of $P$ be denoted by $p_{k\ell}$.
We wish to set multiple nonzero components of $P$ to zero simultaneously.
In particular, let us introduce the set
\begin{align}\label{eq:440}
\mathcal{D}_P\subset\operatorname{supp}(P):=\{(k,\ell)\text{ such that }p_{k\ell}\neq 0\}
\end{align}
of indices corresponding to the entries we wish to drop, where $\operatorname{supp}(P)$ is the set of indices of nonzero components in $P$.
For hard thresholding, whenever $|\operatorname{supp}(P)|>m$ for some $m\ll n^2$, we deploy the following strategy:
\begin{align}\label{eq:445}
\text{For }\;(k,\ell)\in
\mathcal{D}_P=
\arg
\underset{\underset{\mathlarger{|\operatorname{supp}(P)|-|\mathcal{S}|=m}}{\mathcal{S}\subset \operatorname{supp}(P)}}{\min}
\left\{
\sum_{(k,\ell)\in\mathcal{S}}
|p_{k\ell}|
\right\},
\text{ set }\;
p_{k\ell}:=0.
\end{align}

\subsection{Effects of dropping on implementations}
As dropping strategies are applied to the iterates $M_{i+1}$ and $P_i$, the iterative update formulas for the residual iterate $R_{i+1}$ do not hold anymore.
Consequently, for each of the algorithms presented in this paper, once dropping strategies are deployed, the residual iterates need be evaluated explicitly, that is, we use the relation $R_{i+1}:=I_n-AM_{i+1}$ instead of the previously derived update formulas.
In Table~\ref{tab:operation-count}, we detail the operation count per iteration with implicit residual update.
As shown in Table~\ref{tab:operation-count-with-dropping}, when explicit residual evaluation is deployed, each iteration requires an additional product between sparse matrices (SpGEMM).

\begin{table}
\caption{Operation count per iteration without explicit residual update}
\label{tab:operation-count}
\begin{tabular}{rc}
\toprule
Method & Operation count per iteration\\
\midrule
  MR & 1~SpGEMM${}^\dagger$\hphantom{s} + 1~Frobenius inner product\hphantom{s} + 1~Frobenius norm\hphantom{s} + 2~SpGEAMs${}^*$\\
  SD & 2~SpGEMMs\hphantom{${}^\dagger$}\hspace{.02cm} + 1~Frobenius inner product\hphantom{s} + 1~Frobenius norm\hphantom{s} + 2~SpGEAMs\hphantom{${}^*$}\\
  CG & 1~SpGEMM\hphantom{s}\hphantom{${}^\dagger$}\hspace{.02cm} + 1~Frobenius inner product\hphantom{s} + 1~Frobenius norm\hphantom{s} + 3~SpGEAMs\hphantom{${}^*$}\\
LOMR & 2~SpGEMMs\hphantom{${}^\dagger$}\hspace{.02cm} + 3~Frobenius inner products + 2~Frobenius norms + 5~SpGEAMs\hphantom{${}^*$}\\
\bottomrule
\end{tabular}
{\footnotesize$\hspace*{1cm}{}^\dagger$: General product of two sparse matrices.}\\
{\footnotesize$\hspace*{1cm}{}^*$: General addition of two sparse matrices.}
\end{table}

\begin{table}
\caption{Operation count per iteration with explicit residual update}
\label{tab:operation-count-with-dropping}
\begin{tabular}{rc}
\toprule
Method & Operation count per iteration\\
\midrule
  MR & 2~SpGEMM${}^\dagger$\hphantom{s} + 1~Frobenius inner product\hphantom{s} + 1~Frobenius norm\hphantom{s} + 2~SpGEAMs\\
  SD & 3~SpGEMMs\hphantom{${}^\dagger$}\hspace{.02cm} + 1~Frobenius inner product\hphantom{s} + 1~Frobenius norm\hphantom{s} + 2~SpGEAMs\\
  CG & 2~SpGEMM\hphantom{${}^\dagger$}\hspace{.02cm}\hphantom{s} + 1~Frobenius inner product\hphantom{s} + 1~Frobenius norm\hphantom{s} + 3~SpGEAMs\\
LOMR & 3~SpGEMMs\hphantom{${}^\dagger$}\hspace{.02cm} + 3~Frobenius inner products + 2~Frobenius norms + 4~SpGEAMs\\
\bottomrule
\end{tabular}
{\footnotesize$\hspace*{1cm}{}^\dagger$: General product of two sparse matrices.}\\
{\footnotesize$\hspace*{1cm}{}^*$: General addition of two sparse matrices.}
\end{table}

\section{Numerical experiments}\label{sec:numerical-experiments}
\begin{figure}[H]
\centering
\includegraphics{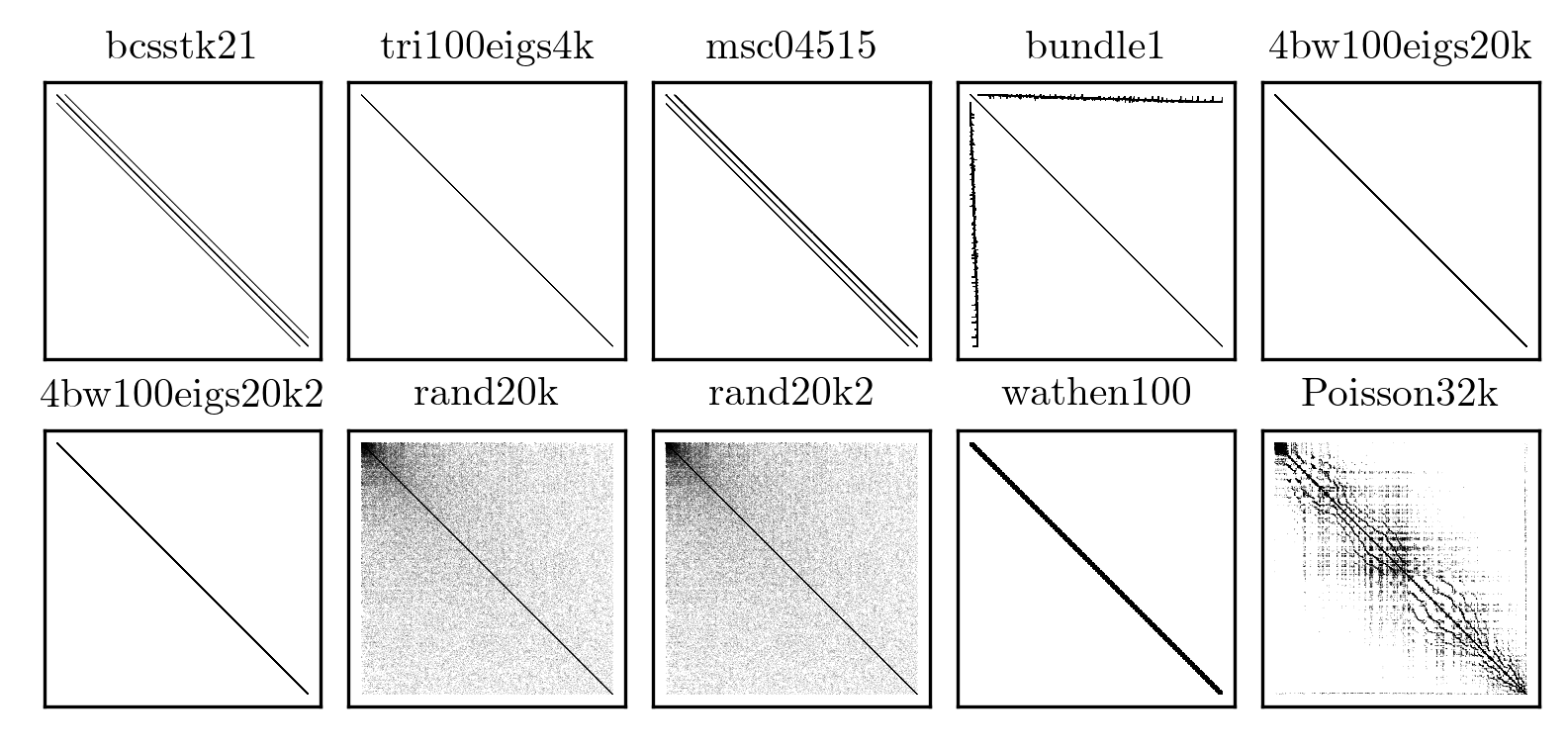}
\caption{Nonzero structures of test sparse matrices (Experiment00).}
\label{fig:sparsity-patterns}
\end{figure}

In this section, we present the results of a number of experiments that were conducted to showcase the behavior of the algorithms presented in this paper, particularly in comparison to the state-of-the-art (P)MR and (P)SD methods that are used for global iterations of approximate inverses. 
These experiments are labeled from Experiment01 to Experiment06 with corresponding Julia scripts and Python post-processing scripts available at 
\begin{center}
\href{https://github.com/venkovic/SparseApproximateInverses.jl/tree/master/examples}{https://github.com/venkovic/SparseApproximateInverses.jl/tree/master/examples}.
\end{center}
For these experiments, we present 10 test sparse SPD matrices and their main characteristics in Table~\ref{tab:metadata-sparse-matrices}.
The nonzero patterns of these matrices are plotted in Fig.~\ref{fig:sparsity-patterns}.
We distinguish between two groups of matrices.
First, we have relatively small matrices, i.e., $n<5{,}000$, for which we are able to compute the entire spectrum, as well as run global iteration algorithms without the need to drop nonzero values.
These matrices are
\begin{itemize}
\item[-] bcsstk21: Ill-conditioned stiffness matrix from a structural mechanics problem of a clamped square plate. 
All eigenvalues are greater than 1.
Bandwidth of 125 with only a handful of nonzero diagonals, typical of finite element 
discretizations.
\item[-] tri100eigs4k: Custom-made ill-conditioned tridiagonal matrix $A=LL^T$, where the lower bidiagonal matrix $L$ has 100 distinct random eigenvalues (diagonal components) independently and identically distributed (iid) as uniform random variables (RVs) in $(0.05,1.05)$. All other eigenvalues (diagonal components) of $L$ are 1. 
The sub-diagonal components of $L$ are iid uniform RVs in $(0,1)$.
$A$ is nearly singular with a smallest eigenvalue of $10^{-8}$.
\item[-] msc04515: Moderately conditioned structural mechanics matrix from Boeing obtained using the NASTRAN finite element software.
Bandwidth of 162 with only a handful of nonzero diagonals.
Relatively large eigenvalues with a spectrum spanning from $10^4$ to $10^{10}$.
\end{itemize}
The seven other matrices, for which the growth of nonzero density needs be contained, have dimensions spanning from more than $10{,}000$ to more than $30{,}000$.
These matrices are
\begin{itemize}
\item[-] bundle1: Moderately conditioned generalized arrowhead matrix from 3D vision bundle adjustment.
Very large eigenvalues spanning from $10^9$ to $10^{12}$.
\item[-] 4bw100eigs20k: Ill-conditioned custom-made banded matrix $A=LL^T$, where the lower triangular matrix $L$ has a dense bandwidth of 4.
The factor $L$ has 100 iid random distinct eigenvalues (diagonal components) uniformly distributed in $(0,0.05)$.
All the other eigenvalues of $L$ are set to 1.
Each sub-diagonal of $L$ has uniform iid random components in $(0,0.01)$.
$A$ is nearly singular with a smallest eigenvalue of $10^{-9}$.
\item[-] 4bw100eigs20k2: Ill-conditioned custom-made matrix $A=LL^T$ similar to 4bw100eigs20k, but for which the random distinct eigenvalues of $L$ are in $(0,10^5)$.
Unlike 4bw100eigs20k, 4bw100eigs20k2 is not nearly singular.
\item[-] rand20k: Ill-conditioned custom-made matrix $A=LL^T$ with random off-diagonal nonzero pattern. 
The lower triangular factor $L$ has iid random eigenvalues in $(0,10^4)$, and all nonzero off-diagonal components are in $(0,1)$.
\item[-] rand20k2: Ill-conditioned custom-made matrix $A=LL^T$ with the same nonzero pattern and off-diagonal values as rand20k.
For rand20k2, the eigenvalues of $L$ are in $(0,10^2)$, rendering $A$ nearly singular, with a smallest eigenvalue of $10^{-6}$.
\item[-] wathen100: Moderately conditioned random banded matrix with a bandwidth of 304, coming from Andy Wathen at Oxford University.
\item[-] Poisson32k:
Moderately conditioned stiffness matrix of a P1 finite element unstructured discretization of the 2D Poisson PDE with a random variable coefficient sampled as a Gaussian process with unit-variance squared exponential covariance.
Originally from~\cite{venkovic2023preconditioning}.
\end{itemize}

\begin{table}
\caption{Metadata of test sparse matrices}
\label{tab:metadata-sparse-matrices}
\begin{tabular}{ccccccc}
\toprule
Matrix & $n$ & $nnz$ & $nnz/n^2$ & $\lambda_{min}$ & $\lambda_{max}$ & $\kappa$\\
\midrule
bcsstk21${}^{\dagger}$ & 3,600 & 26,600 & $2.05\times 10^{-3}$ & $7.21\times 10^{0\hphantom{-}}$ & $1.27\times 10^{8\hphantom{0}}$ & $1.76\times 10^{7\hphantom{0}}$ \\
tri100eigs4k${}^{*}$ & 4,000 & 11,998  & $7.50\times 10^{-4}$ & $9.26\times 10^{-9}$ & $3.56\times 10^{0\hphantom{0}}$&  $3.85\times 10^{8\hphantom{0}}$ \\
msc04515${}^{*}$ & 4,515 & 97,707 & $4.79\times 10^{-3}$ & $1.39\times 10^{4\hphantom{-}}$ & $3.15\times 10^{10}$ & $2.26\times 10^{6\hphantom{0}}$\\
\midrule
bundle1${}^\dagger$ & 10,581 & 770,811 & $6.89\times 10^{-3}$ & $6.40\times 10^{9\hphantom{-}}$ & $6.43\times 10^{12}$ & $1.00\times 10^{3\hphantom{0}}$ \\
4bw100eigs20k${}^{*}$ & 20,000 & 179,980 & $4.50\times 10^{-4}$  & $2.47\times 10^{-9}$ & $1.05\times 10^{0\hphantom{0}}$ & $4.25\times 10^{8\hphantom{0}}$ \\
4bw100eigs20k2${}^{*}$ & 20,000 & 179,980 & $4.50\times 10^{-4}$ & $9.74\times 10^{-1}$ & $9.79\times 10^{9\hphantom{0}}$ & $1.01\times 10^{10}$ \\
rand20k${}^{*}$ & 20,000 & 99,772 & $2.49\times 10^{-4}$&$8.70\times 10^{-2}$&$1.00\times 10^{8\hphantom{0}}$&$1.15\times 10^{9\hphantom{0}}$\\
rand20k2${}^{*}$ & 20,000 & 99,772 & $2.49\times 10^{-4}$ & $8.67\times 10^{-6}$ & $1.00\times 10^{4\hphantom{0}}$&$1.15\times 10^{9\hphantom{0}}$\\	
wathen100${}^\dagger$ & 30,401 & 471,601 & $5.10\times 10^{-4}$ & $6.36\times 10^{-2}$ & $3.70\times 10^{2\hphantom{0}}$ & $5.82\times 10^{3\hphantom{0}}$ \\
Poisson32k${}^*$ & 31,839 & 221,375 & $2.18\times 10^{-4}$ & $6.21\times 10^{-4}$ & $9.66\times 10^{2\hphantom{0}}$ & $1.55\times 10^{5\hphantom{0}}$ \\
\bottomrule
\end{tabular}
{\footnotesize$\hspace*{1cm}{}^\dagger$: \url{https://sparse.tamu.edu/}}\\
{\footnotesize$\hspace*{1cm}{}^*$: \url{https://github.com/venkovic/matrix-market}}
\end{table}

\subsection{Dropping-free experiments}\label{sec:dropping-free-experiments}
In this section, we present the results of dropping-free experiments, first on small matrices, for which entire spectra can be computed, and then, on medium size matrices.

\subsubsection{Small matrices}\label{sec:small-matrices}
We benchmark four methods methods (PSD\_SPD, PMR\_SPD, PCG, and LOPMR\_SPD) using Jacobi preconditioners to compute approximate inverses of three small matrices: bcsstk21, tri100eigs4k, and msc04515. 
This experiment (Experiment01 in the corresponding GitHub repository) is performed without dropping, allowing the nonzero densities of the main iterate and search direction to grow uncontrollably.
Despite this unconstrained expansion of the nonzero pattern, the approximate inverse of each matrix achieve densities that are  actually bounded in relation to the sparsity pattern of $A$, the initial matrix iterate, and the nature of the global iteration method. 
For bcsstk21, all four methods yield approximate inverses with densities of exactly 50\%. 
For msc04515, all approximate inverses achieve 52\% density. 
For the tri100eigs4k matrix, PMR\_SPD, LOPMR\_SPD and PCG all achieve a density of 14.5\%, whereas PSD\_SPD yields a density of 22.3\%.

\begin{figure}[H]
\centering
\includegraphics[scale=.8]{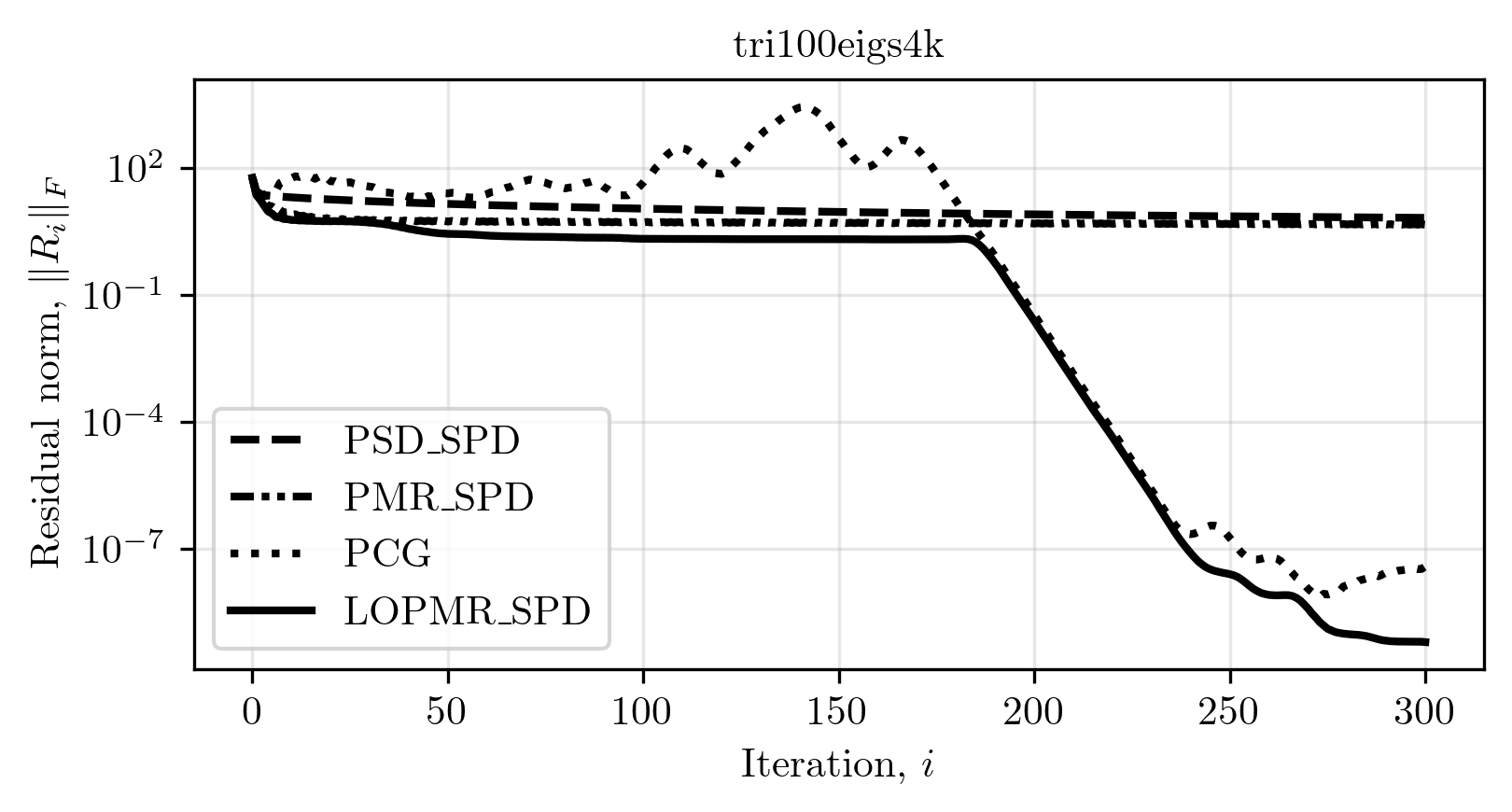}
\includegraphics[scale=.8]{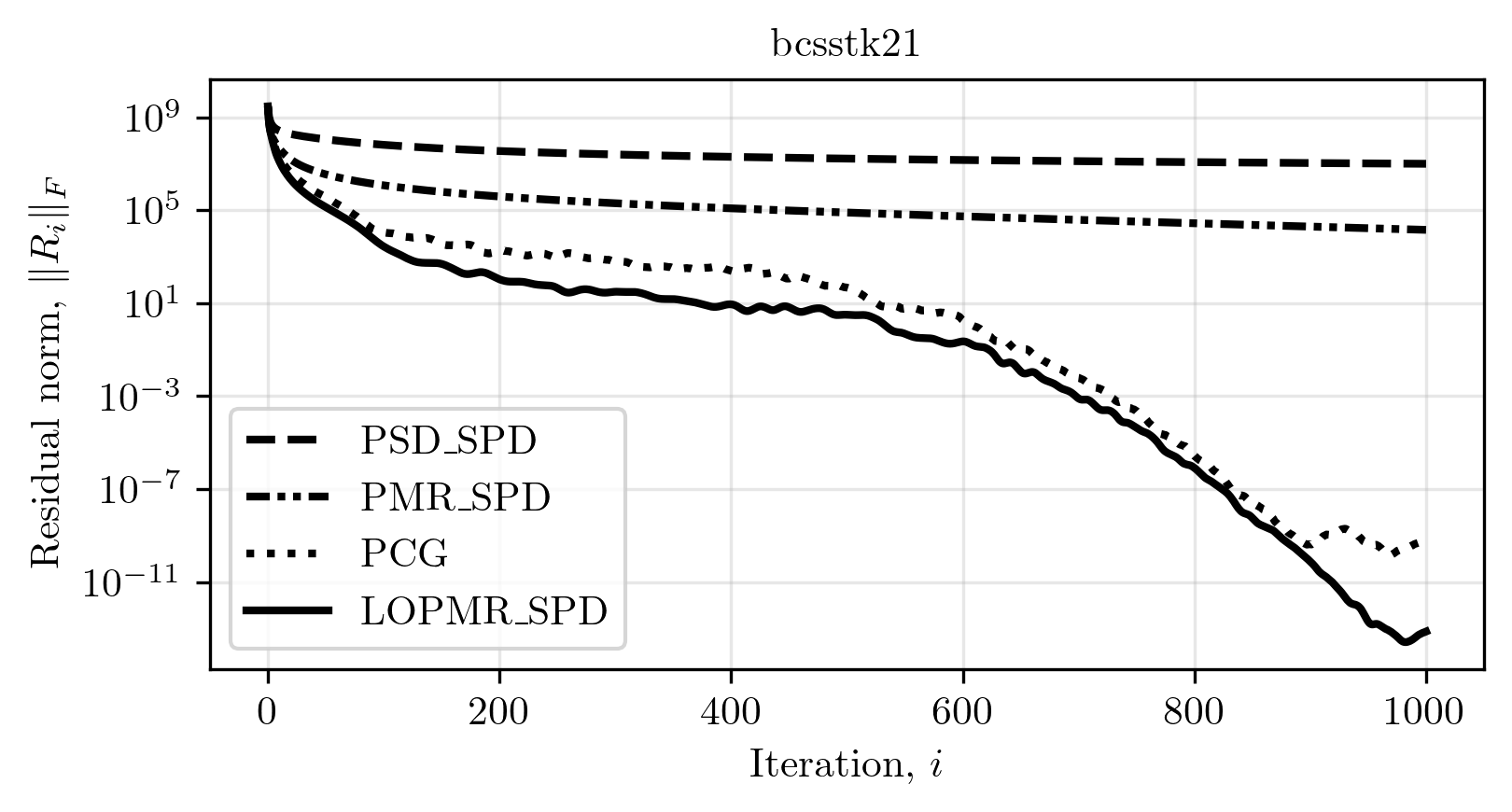}
\includegraphics[scale=.8]{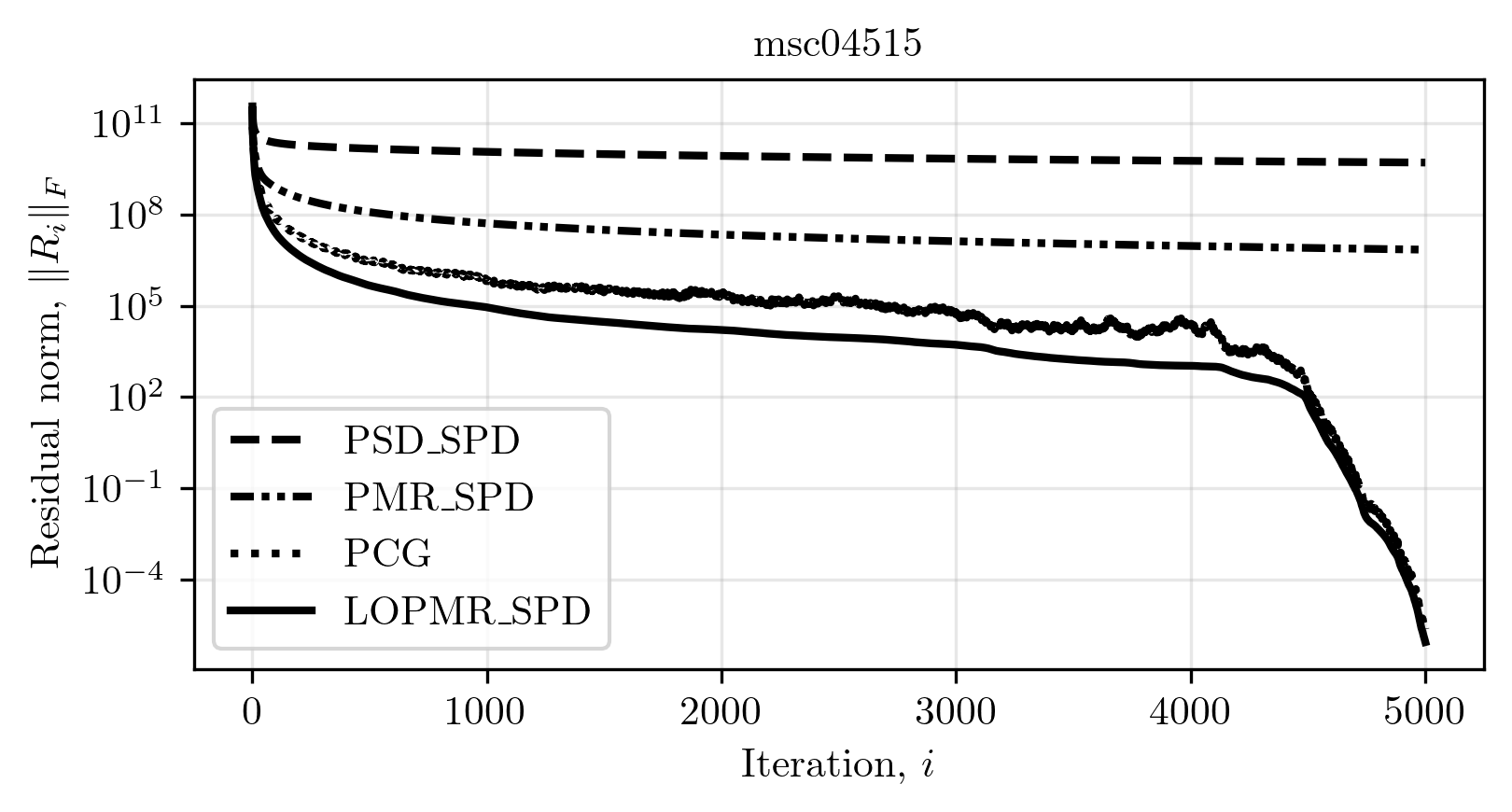}
\caption{Convergence plots of dropping-free experiments (Experiment01).}
\label{fig:convergence-dropping-free}
\end{figure}

We plot the Frobenius norm of the residual as a function of the number of matrix iterations in Figure~\ref{fig:convergence-dropping-free}. 
The results show that both PCG and LOPMR\_SPD achieve convergence for all three tested matrices. 
LOPMR\_SPD exhibits monotone decrease of the residual norm, whereas PCG, which does not minimize the residual norm, occasionally shows significant oscillations which are particularly evident for the nearly singular and ill-conditioned matrix tri100eigs4k.
Note, however, that after several hundred iterations on the tri100eigs4k matrix, PCG does converge and nearly matches LOPMR\_SPD's behavior. 
This convergence pattern is consistent with Krylov subspace approximations of matrices having few distinct eigenvalues, which possess low-degree minimal polynomials that enable rapid convergence once the iteration count approaches the number of distinct eigenvalues.

Figure~\ref{fig:convergence-dropping-free} shows that the Frobenius residual norm barely converges for PSD\_SPD and PMR\_SPD. 
Among these two methods, PSD\_SPD performs worst. 
In all cases, these two methods fail to achieve the symbolic threshold of $\|R\|_F<1$, below which the approximate inverse is guaranteed to be non-singular (though not necessarily SPD).

\begin{figure}[H]
\centering
\includegraphics[scale=.8]{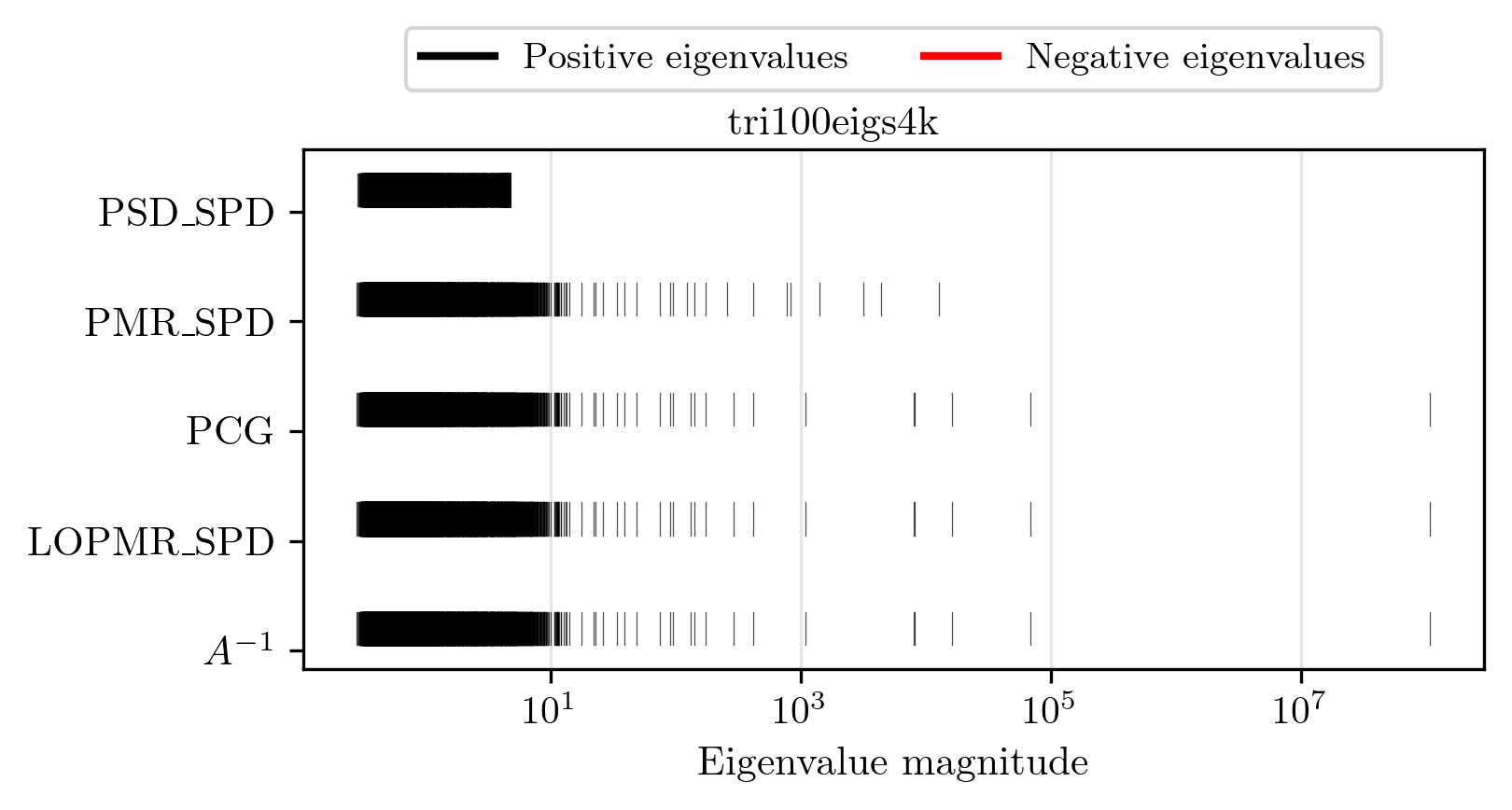}
\includegraphics[scale=.8]{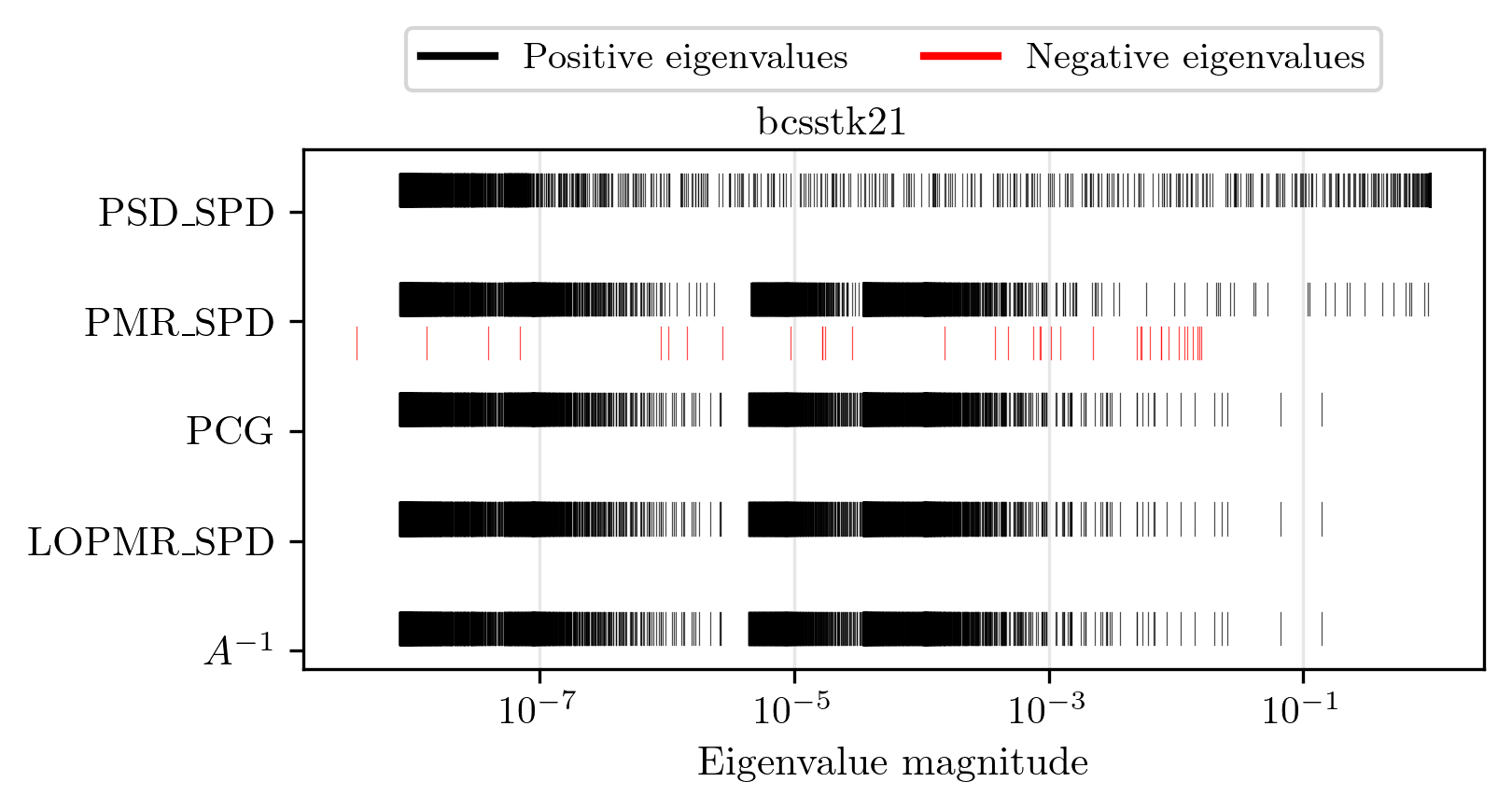}
\includegraphics[scale=.8]{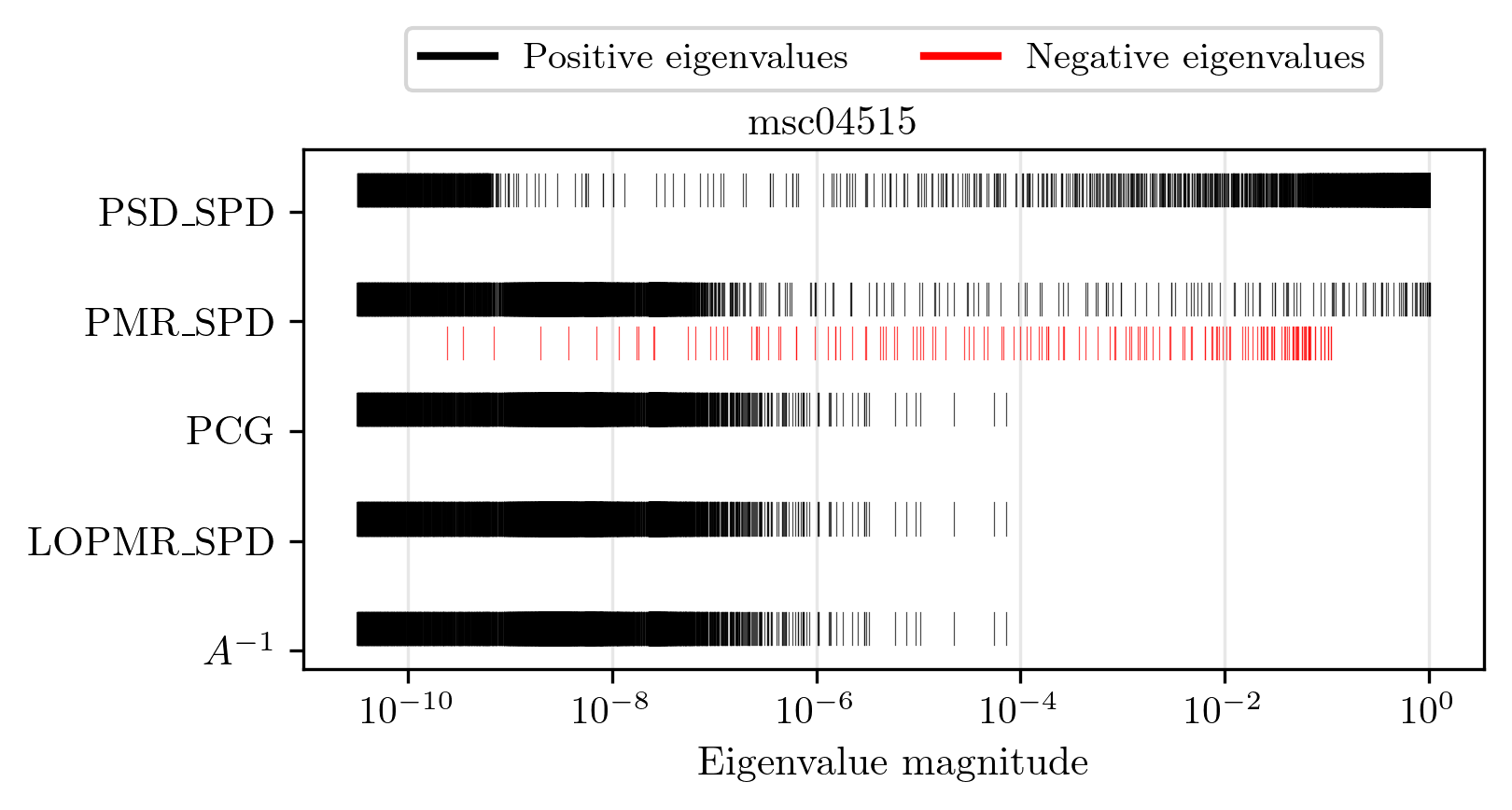}
\caption{Converged spectra of dropping-free experiments (Experiment01).}
\label{fig:spectra-dropping-free}
\end{figure}

Figure~\ref{fig:spectra-dropping-free} presents the spectra of the approximate inverses obtained by each method for all three small matrices. 
Most notably, the converged approximate inverses of PCG and LOPMR\_SPD remarkably reproduce the spectrum of $A^{-1}$. 
Importantly, these approximate inverses are SPD, which has significant implications for deploying an approximate inverse as a preconditioner in standard vector CG iterations.
In contrast, the other two methods, namely PSD\_SPD and PMR\_SPD, both fail to reproduce the spectrum of $A^{-1}$.
In particular, PMR\_SPD fails to produce positive definite approximate inverses, which represents a major shortcoming for practical applications. 
Note, however, that we purposely selected matrices that expose the limitations of these methods. 
For some matrices, all four methods can produce SPD approximate inverses that reasonably reproduce the spectrum of $A^{-1}$.
We invite readers to test this with the Poisson4k\footnotemark[1], triclust4k\footnotemark[1] and triunif4k\footnotemark[1] matrices. 
Nevertheless, both PCG and LOPMR\_SPD consistently outperform the other two methods.
\footnotetext[1]{\url{https://github.com/venkovic/matrix-market}}

\subsubsection{Medium size matrices}\label{sec:medium-size-matrices}
We now focus on the seven other matrices for which we cannot afford to let the approximate inverse and search direction grow uncontrollably dense, nor can we compute complete spectra for these matrices. 
Nevertheless, we first examine how global iteration methods perform without dropping. 
Given the behaviors observed in the small matrix experiments, we hereafter focus on the two best-performing methods (LOPMR\_SPD and PCG) as well as the current state-of-the-art method for global iteration (PMR\_SPD).
We deploy these three methods with a stopping criterion based on the density of the main iterate $M_i$: iteration stops once a density of 3\% or higher is achieved. 
Given the different nonzero patterns of the matrices under consideration, the number of iterations required to reach this stopping criterion varies significantly, as reported in Table~\ref{tab:dropping-free-summary}.

Table~\ref{tab:dropping-free-summary} shows that, besides the custom-made matrix 4bw100eigs20k2, all methods produce iterates with exactly the same density. 
Another notable observation is how far beyond the 3\% threshold some iterates can jump in a single iteration. 
For instance, the approximate inverses for bundle1 all achieve 21.5\% density. 
Less dramatically but still significantly, the approximate inverses of rand20k and rand20k2 reach 7.86\% density.
We attribute bundle1's sudden density increase to the matrix's arrowhead structure, while for the two random matrices, we assume the randomness of the nonzero pattern causes the abrupt rise in main iterate density. 
These observations are important because they indicate which matrices will likely require more aggressive dropping in the next section. 
We will see that large amounts of dropping not only carry computational cost but can also detrimentally impact the convergence of global iteration methods.

\begin{figure}[ht]
\centering
\includegraphics[scale=.8]{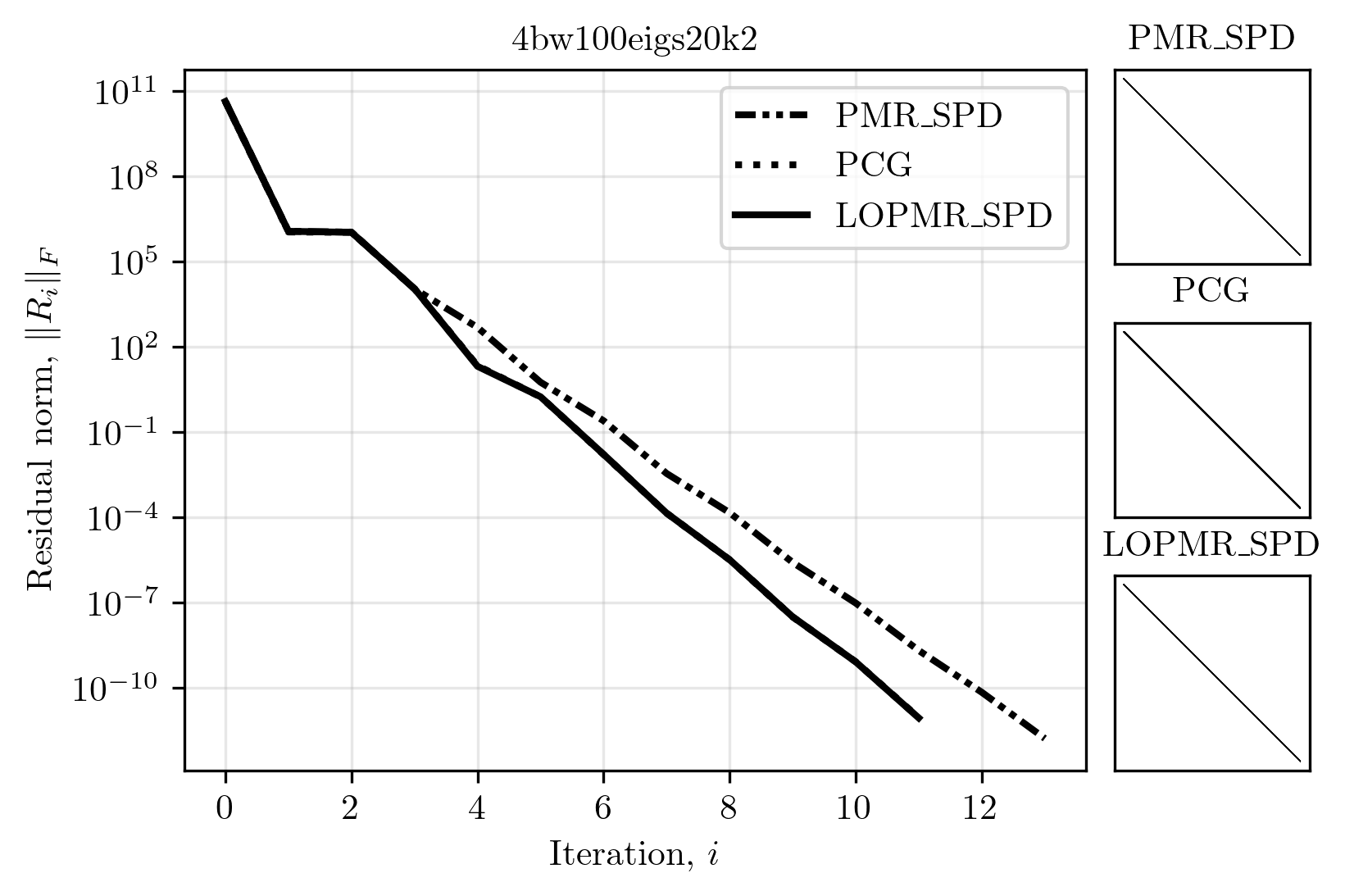}
\includegraphics[scale=.8]{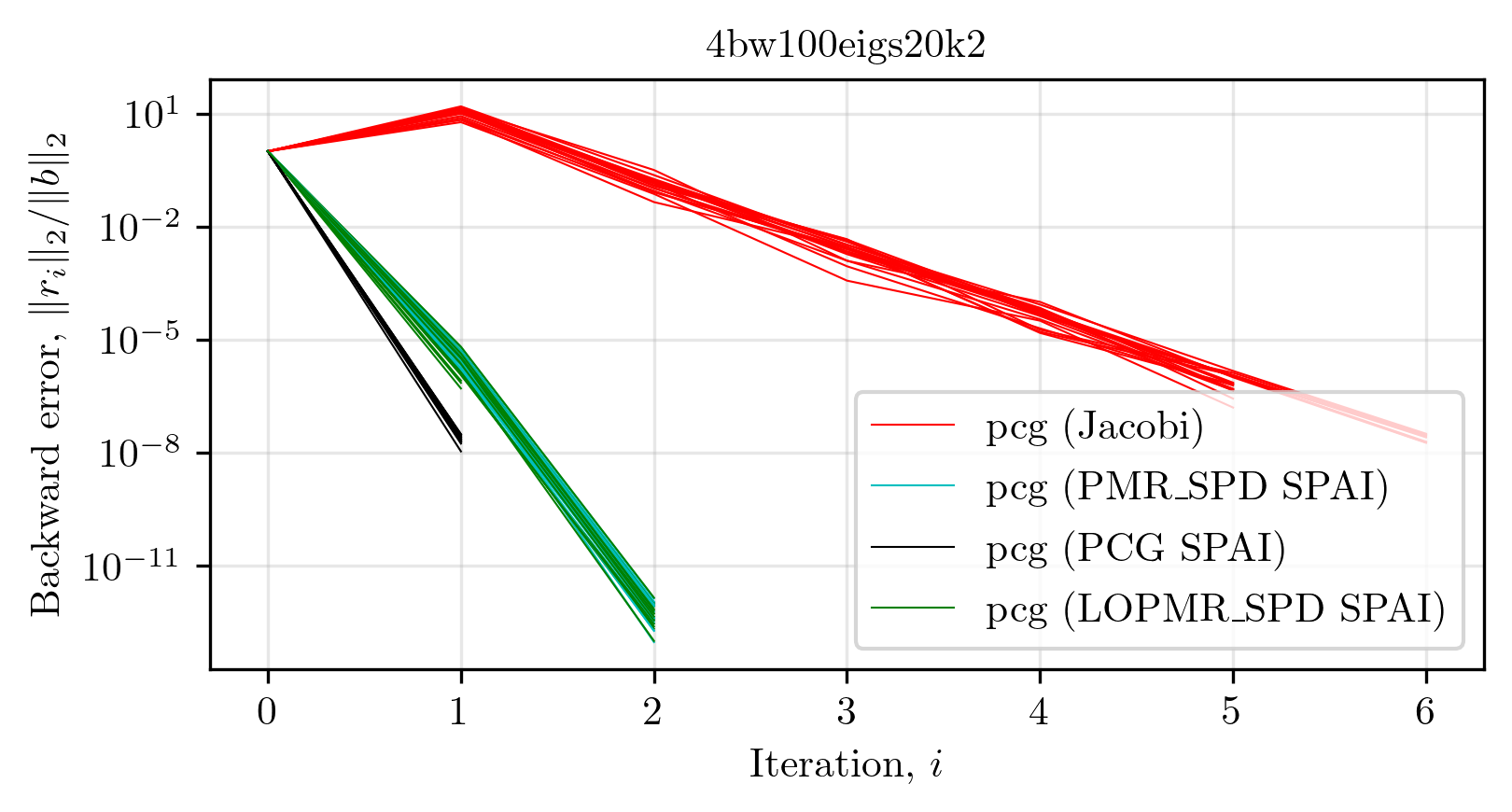}
\caption{Convergence plot of dropping-free SPAIs and PCG results for the matrix 4bw100eigs20k2 (Experiment03-04).}
\label{fig:4bw100eigs20k2}
\end{figure}

Except for 4bw100eigs20k2, all tested matrices reach the stopping criterion based on the density upper bound before converging to a proper approximate inverse. 
Despite having the worst condition number, the Frobenius residual norm of the 4bw100eigs20k2 approximate inverse decreases very quickly for all three methods, see top insert in Figure~\ref{fig:4bw100eigs20k2}, achieving $10^{-6}$ in no more than ten matrix iterations irrespective of the method used.
Consequently, all three methods stop well before reaching the 3\% density threshold, with a density of $0.524\%$ for PMR\_SPD and $0.444\%$ for PCG and LOPMR\_SPD.
Given this behavior, 4bw100eigs20k2 is not a suitable candidate for dropping experiments.
The quality of the approximate inverses generated for 4bw100eigs20k2 is assessed in the bottom part of Figure~\ref{fig:4bw100eigs20k2}, where the convergence of vector PCG iterations to a $10^{-6}$ backward error is shown using these approximate inverses as preconditioners in comparison to using a Jacobi preconditioner. 
The approximate inverses generated by PMR\_SPD, PCG and LOPMR\_SPD  all constitute excellent preconditioners, reducing convergence from 6 to 1 or 2 iterations, with the same densities around $0.5\%$. 

For the six other medium-sized matrices, most approximate inverses are indefinite, and those that are positive definite do not yet constitute good preconditioners. 
Therefore, all six matrices (bundle1, 4bw100eigs20k, rand20k, rand20k2, wathen100, and Poisson32k) are suitable candidates for dropping experiments. 
By conducting these experiments, we aim to continue global iterations to further decrease the residual norm while dropping sufficient nonzero values to control memory usage. 
These experiments are documented in Section~\ref{sec:dropping-experiments}.

\begin{table}
\caption{Summary results of dropping-free SPAIs stopped after density reaches 3\% (Experiment03).}
\label{tab:dropping-free-summary}
\begin{tabular}{cccccc}
\toprule
Method & $\lambda_{min}$                  & $\lambda_{max}$     & $\|I_n-AM\|_F$ 
& \# $iters$ & $nnz/n^2$\\
\midrule
\multicolumn{6}{c}{bundle1}\\
PMR\_SPD 
& $-6.96\times 10^{0}$ & $8.24\times 10^{0}$ & $4.36\times 10^{12}$ & 2 & $2.15\times 10^{-1}$\\
PCG
& $-7.39\times 10^{0}$ & $8.29\times 10^{0}$ & $3.97\times 10^{12}$ & 2 & $2.15\times 10^{-1}$\\
LOPMR\_SPD
& $-7.79\times 10^{0}$ & $8.82\times 10^{0}$ & $2.80\times 10^{12}$ & 2 & $2.15\times 10^{-1}$\\
\midrule
\multicolumn{6}{c}{4bw100eigs20k}\\
PMR\_SPD    & $\hphantom{-}1.92\times 10^{0}$  & $3.94\times 10^{6}$ & $3.09\times 10^0$    
& 76   & $3.02\times 10^{-2}$\\
PCG    & $-2.60\times 10^{5}$             & $6.10\times 10^{8}$ & $3.61\times 10^{2}$  
& 76 & $3.02\times 10^{-2}$\\
LOPMR\_SPD  & $\hphantom{-}1.91\times 10^{0}$  & $7.10\times 10^{7}$ & $1.41\times 10^0$    
& 76 & $3.02\times 10^{-2}$\\
\midrule
\multicolumn{6}{c}{4bw100eigs20k2}\\
PMR\_SPD & $\hphantom{-}5.06\times 10^{-9}$ & $2.05\times 10^{0}$ & $1.68\times 10^{-12}$ 
& 13  & $5.24\times 10^{-3}$\\
PCG    & $\hphantom{-}3.82\times 10^{-9}$ & $2.05\times 10^{0}$ & $9.11\times 10^{-12}$ 
& 11    & $4.44\times 10^{-3}$\\
LOPMR\_SPD  & $\hphantom{-}3.07\times 10^{-9}$ & $2.05\times 10^{0}$ &$8.87\times 10^{-12}$ 
& 11    & $4.44\times 10^{-3}$\\
\midrule
\multicolumn{6}{c}{rand20k}\\
PMR\_SPD & $-2.17\times 10^{2}$           & $2.22\times 10^{2}$ & $1.77\times 10^5$ 
& 5    & $7.86\times 10^{-2}$\\
PCG    & $-3.10\times 10^{1}$             & $3.59\times 10^{1}$ & $8.07\times 10^4$ 
& 5    & $7.86\times 10^{-2}$\\
LOPMR\_SPD  & $-4.97\times 10^{1}$             & $5.46\times 10^{1}$ & $1.19\times 10^5$ 
& 5    & $7.86\times 10^{-2}$\\
\midrule
\multicolumn{6}{c}{rand20k2}\\
PMR\_SPD & $-9.18\times 10^{2}$             & $8.81\times 10^{3}$ & $1.74\times 10^3$ 
& 5    & $7.86\times 10^{-2}$\\
PCG    & $-8.37\times 10^{1}$             & $8.81\times 10^{3}$ & $1.94\times 10^3$ 
& 5    & $7.86\times 10^{-2}$\\
LOPMR\_SPD  & $-8.68\times 10^{1}$             & $8.81\times 10^{3}$ & $1.23\times 10^3$ 
& 5    & $7.86\times 10^{-2}$\\
\midrule
\multicolumn{6}{c}{wathen100}\\
PMR\_SPD    & $-1.36\times 10^{-1}$            & $3.11\times 10^{1}$ & $7.18\times 10^2$ 
& 9    & $3.02\times 10^{-2}$\\
PCG    & $-4.26\times 10^{-2}$            & $3.11\times 10^{1}$ & $1.26\times 10^2$ 
& 9    & $3.02\times 10^{-2}$\\
LOPMR\_SPD  & $-1.28\times 10^{-2}$            & $3.10\times 10^{1}$ & $1.04\times 10^2$ 
& 9    & $3.02\times 10^{-2}$\\
\midrule
\multicolumn{6}{c}{Poisson32k}\\
PMR\_SPD    & $-6.56\times 10^{-3}$            & $1.03\times 10^{2}$ & $1.48\times 10^1$ 
& 17   & $3.22\times 10^{-2}$\\
PCG    & $-8.20\times 10^{-1}$            & $5.06\times 10^{2}$ & $3.35\times 10^{1}$ 
& 17   & $3.22\times 10^{-2}$\\
LOPMR\_SPD  & $\hphantom{-}1.96\times 10^{-2}$ & $3.44\times 10^{2}$ & $7.99\times 10^0$ 
& 17   & $3.22\times 10^{-2}$\\
\bottomrule
\end{tabular}
\end{table}

As previously mentioned, our sparse matrix selection is designed to demonstrate situations where PCG and, more particularly, LOPMR\_SPD outperform PMR\_SPD. 
However, some matrices allow all three methods to construct SPD approximate inverses that constitute good preconditioners even when limiting the main iterate density to a maximum of 3\%. 
This is the case for the following matrices:
\begin{itemize}
\item[-] crystm02\footnotemark[2]: Well-conditioned banded mass matrix for the finite element modeling of free crystal vibration.
Submitted by Boeing.
Nearly singular, with a smallest eigenvalue of $10^{-15}$.
\item[-] Dubcova1\footnotemark[2]: Well-conditioned matrix obtained with a PDE solver.
Submitted by Dubcova et al. from the University of Texas at El Paso.
\item[-] jnlbrng1\footnotemark[2]: Well-conditioned banded matrix from a quadratic journal bearing optimization problem.
\item[-] minsurfo\footnotemark[2]: Well-conditioned banded matrix from minimum surface optimization problem.
\item[-] ted\_B\footnotemark[2]: Ill-conditioned banded matrix obtained by finite element method applied to coupled thermoelasticity equations on a bar.
Nearly singular matrix, with a smallest eigenvalue of $10^{-8}$.
Submitted by Bindel from UC Berkeley.
\end{itemize}
Note that for all these matrices, the approximate inverse obtained by LOPMR\_SPD always performed at least as a well as a preconditioner to vector PCG iterations than the approximate inverses obtained by PMR\_SPD and PCG.
\footnotetext[2]{\url{https://sparse.tamu.edu/}}

\subsection{Experiments with dropping}\label{sec:dropping-experiments}
Experiment05 through Experiment07 are performed with nonzero dropping strategies.
In Experiment05, we apply the PMR\_SPD, PCG, and LOPMR\_SPD methods to compute SPAIs using the dropping strategy described in Section~\ref{sec:dropping}. 
The results are summarized in Tables~\ref{tab:with-dropping-summary-1}-\ref{tab:with-dropping-summary-2} after 5, 20 and 100 matrix iterations, unless convergence in terms of backward error of $10^{-6}$ is achieved faster.
From the results of Experiment05, we note that dropping takes effect from the first iteration through symmetrization of the main iterate, followed by dropping nonzero values with magnitudes smaller than unit round-off.
Consequently, even though the 4bw100eigs20k and rand20k matrices quickly reached the maximum 3\% density in the previous experiment, they do not reach this maximum density in Experiment05. 
From Tables~\ref{tab:with-dropping-summary-1}-\ref{tab:with-dropping-summary-2}, we see that we are able to build SPAIs with 3\% density and backward errors of $10^{-6}$ for half of the listed matrices, namely 4bw100eigs20k, rand20k and wathen100, using both PCG and LOPMR\_SPD and, for rand20k, also when using PMR\_SPD.
In Experiment06, we assess the performance of the SPAIs obtained from Experiment05 when used as preconditioners for linear solves by CG vector iteration.
While achieving convergence of the matrix iteration is one way to build an SPAI that performs well as a preconditioner, it is nevertheless possible to achieve good preconditioning with an SPAI even if it originates from a non-converged matrix iteration, see, e.g., the case of Poisson32k.
From here, we go one-by-one through the six matrices listed in Tables~\ref{tab:with-dropping-summary-1}-\ref{tab:with-dropping-summary-2} and comment on the outcome of Experiment05 and Experiment06.

\begin{table}
\caption{Summary results of SPAIs obtained with dropping of nonzero values for maximum densities of 3\% (Experiment05).}
\label{tab:with-dropping-summary-1}
\begin{tabular}{cccccc}
\toprule
Method & $\lambda_{min}$                   & $\lambda_{max}$                  & $\|I_n-AM\|_F$             &
\# $iters$ & $nnz/n^2$\\
\midrule
\multicolumn{6}{c}{bundle1}\\
\multirow{3}{*}{PMR\_SPD} 
& $-5.74\times 10^{-2}$ & $7.20\times 10^{-1}$ & $1.25\times 10^{12}$ & 5 & $3.00\times 10^{-2}$\\
& $-8.02\times 10^{-4}$ & $2.96\times 10^{-1}$ & $1.48\times 10^{11}$ & 20 & $3.00\times 10^{-2}$\\
& $-1.42\times 10^{-6}$ & $5.41\times 10^{-3}$ & $2.02\times 10^{9\hphantom{1}}$ & 100 & $3.00\times 10^{-2}$\\
\cmidrule{2-6}
\multirow{3}{*}{PCG} 
& $-7.31\times 10^{2}$ & $5.65\times 10^{-1}$ & $6.52\times 10^{11}$ & 5 & $3.00\times 10^{-2}$\\
& $-9.63\times 10^{-3}$ & $2.63\times 10^{-3}$ & $6.52\times 10^{9\hphantom{1}}$ & 20 & $3.00\times 10^{-2}$\\
& $-4.23\times 10^{-7}$ & $3.45\times 10^{-7}$ & $4.13\times 10^{5\hphantom{1}}$ & 100 & $3.00\times 10^{-2}$\\
\cmidrule{2-6}
\multirow{3}{*}{LOPMR\_SPD} 
& $-5.04\times 10^{-2\hphantom{1}}$ & $5.98\times 10^{-1\hphantom{1}}$ & $4.76\times 10^{11}$ & 5 & $3.00\times 10^{-2}$\\
& $-1.21\times 10^{-2\hphantom{1}}$ & $1.19\times 10^{-2\hphantom{0}}$ & $4.43\times 10^{9\hphantom{1}}$ & 20 & $3.00\times 10^{-2}$\\
& $\hphantom{-}1.67\times 10^{-13}$ & $1.57\times 10^{-10}$ & $1.96\times 10^{1\hphantom{1}}$ & 100 & $3.00\times 10^{-2}$\\
\midrule
\multicolumn{6}{c}{4bw100eigs20k}\\
\multirow{3}{*}{PMR\_SPD} 
& $\hphantom{-}9.60\times 10^{-1}$ & $2.36\times 10^{5}$ & $5.38\times 10^{0\hphantom{-}}$ & 5 & $2.05\times 10^{-3}$\\
& $\hphantom{-}9.59\times 10^{-1}$ & $6.16\times 10^{5}$ & $3.86\times 10^{0\hphantom{-}}$ & 20 & $2.65\times 10^{-3}$\\
& $\hphantom{-}9.59\times 10^{-1}$ & $2.52\times 10^{6}$ & $2.99\times 10^{0\hphantom{-}}$ & 100 & $2.66\times 10^{-3}$\\
\cmidrule{2-6}
\multirow{3}{*}{PCG} 
& $-2.06\times 10^{3\hphantom{-}}$ & $1.25\times 10^{6}$ & $2.64\times 10^{1\hphantom{-}}$ & 5 & $2.05\times 10^{-3}$\\
& $-4.37\times 10^{4\hphantom{-}}$ & $1.65\times 10^{7}$ & $2.38\times 10^{1\hphantom{-}}$ & 20 & $3.48\times 10^{-3}$\\
& $\hphantom{-}9.59\times 10^{-1}$ & $4.04\times 10^{8}$ & $1.29\times 10^{-4}$ & 97 & $4.55\times 10^{-3}$\\
\cmidrule{2-6}
\multirow{3}{*}{LOPMR\_SPD} 
& $\hphantom{-}9.64\times 10^{-1}$ & $3.87\times 10^{5\hphantom{1}}$ & $4.57\times 10^{0\hphantom{-}}$ & 5 & $2.05\times 10^{-3}$\\
& $\hphantom{-}9.61\times 10^{-1}$ & $5.58\times 10^{6\hphantom{0}}$ & $3.30\times 10^{0\hphantom{-}}$ & 20 & $3.31\times 10^{-3}$\\
& $\hphantom{-}9.59\times 10^{-1}$ & $4.04\times 10^{8}$ & $1.39\times 10^{-4}$ & 97 & $4.55\times 10^{-3}$\\
\midrule
\multicolumn{6}{c}{rand20k}\\
\multirow{3}{*}{PMR\_SPD} 
& $\hphantom{-}1.22\times 10^{-6}$ & $1.08\times 10^{1}$ & $1.77\times 10^{5\hphantom{-}}$ & 5 & $1.55\times 10^{-3}$\\
& $\hphantom{-}1.19\times 10^{-6}$ & $1.15\times 10^{1}$ & $6.90\times 10^{1\hphantom{-}}$ & 20 & $2.51\times 10^{-4}$\\
& $\hphantom{-}1.21\times 10^{-6}$ & $1.15\times 10^{1}$ & $1.13\times 10^{-4}$ & 41 & $2.52\times 10^{-4}$\\
\cmidrule{2-6}
\multirow{3}{*}{PCG} 
& $-1.15\times 10^{-3}$ & $1.13\times 10^{1}$ & $8.07\times 10^{4\hphantom{-}}$ & 5 & $1.75\times 10^{-3}$\\
& $\hphantom{-}1.18\times 10^{-6}$ & $1.14\times 10^{1}$ & $4.88\times 10^{-6}$ & 17 & $2.60\times 10^{-4}$\\
\cmidrule{2-6}
\multirow{3}{*}{LOPMR\_SPD} 
& $-1.40\times 10^{-3}$ & $1.12\times 10^{1}$ & $1.19\times 10^{5\hphantom{-}}$ & 5 & $1.77\times 10^{-3}$\\
& $\hphantom{-}1.29\times 10^{-6}$ & $1.15\times 10^{1}$ & $7.39\times 10^{-6}$ & 17 & $2.54\times 10^{-4}$\\
\bottomrule
\end{tabular}
\end{table}

\begin{table}
\caption{Summary results of SPAIs obtained with dropping of nonzero values for maximum densities of 3\% (Experiment05).}
\label{tab:with-dropping-summary-2}
\begin{tabular}{cccccc}
\toprule
Method & $\lambda_{min}$                   & $\lambda_{max}$                  & $\|I_n-AM\|_F$             &
\# $iters$ & $nnz/n^2$\\
\midrule
\multicolumn{6}{c}{rand20k2}\\
\multirow{3}{*}{PMR\_SPD} 
& $-5.29\times 10^{-2}$ & $4.41\times 10^{3}$ & $1.74\times 10^{3}$ & 5 & $3.00\times 10^{-2}$\\
& $\hphantom{-}3.66\times 10^{-3}$ & $4.41\times 10^{3}$ & $2.48\times 10^{2}$ & 20 & $2.04\times 10^{-2}$\\
& $\hphantom{-}3.99\times 10^{-3}$ & $4.41\times 10^{3}$ & $1.26\times 10^{2}$ & 100 & $3.85\times 10^{-3}$\\
\cmidrule{2-6}
\multirow{3}{*}{PCG} 
& $-7.27\times 10^{-1}$ & $4.40\times 10^{3}$ & $1.94\times 10^{3}$ & 5 & $3.00\times 10^{-2}$\\
& $-1.68\times 10^{-1}$ & $4.41\times 10^{3}$ & $1.61\times 10^{3}$ & 20 & $3.00\times 10^{-2}$\\
& $-6.96\times 10^{1\hphantom{-}}$ & $9.64\times 10^{4}$ & $7.36\times 10^{2}$ & 100 & $3.00\times 10^{-2}$\\
\cmidrule{2-6}
\multirow{3}{*}{LOPMR\_SPD} 
& $\hphantom{-}5.52\times 10^{-3}$ & $4.41\times 10^{3}$ & $1.23\times 10^{3}$ & 5 & $3.00\times 10^{-2}$\\
& $\hphantom{-}3.88\times 10^{-3}$ & $4.41\times 10^{3}$ & $1.35\times 10^{2}$ & 20 & $3.00\times 10^{-2}$\\
& $\hphantom{-}4.05\times 10^{-3}$ & $5.30\times 10^{4}$ & $6.41\times 10^{1}$ & 100 & $3.00\times 10^{-2}$\\
\midrule
\multicolumn{6}{c}{wathen100}\\
\multirow{3}{*}{PMR\_SPD} 
& $-9.83\times 10^{-2}$ & $1.43\times 10^{1}$ & $1.55\times 10^{3\hphantom{-}}$ & 5 & $9.99\times 10^{-3}$\\
& $\hphantom{-}3.40\times 10^{-3}$ & $1.57\times 10^{1}$ & $1.35\times 10^{2\hphantom{-}}$ & 20 & $3.00\times 10^{-2}$\\
& $\hphantom{-}2.71\times 10^{-3}$ & $1.57\times 10^{1}$ & $6.98\times 10^{-3}$ & 100 & $3.00\times 10^{-2}$\\
\cmidrule{2-6}
\multirow{3}{*}{PCG} 
& $-1.23\times 10^{-1}$ & $1.71\times 10^{1}$ & $8.19\times 10^{2\hphantom{-}}$ & 5 & $9.99\times 10^{-3}$\\
& $\hphantom{-}2.70\times 10^{-3}$ & $1.57\times 10^{1}$ & $6.64\times 10^{-1}$ & 20 & $3.00\times 10^{-2}$\\
& $\hphantom{-}2.71\times 10^{-3}$ & $1.57\times 10^{1}$ & $1.15\times 10^{-4}$ & 40 & $3.00\times 10^{-2}$\\
\cmidrule{2-6}
\multirow{3}{*}{LOPMR\_SPD} 
& $-8.16\times 10^{-2}$ & $1.71\times 10^{1}$ & $6.73\times 10^{2\hphantom{-}}$ & 5 & $1.00\times 10^{-2}$\\
& $\hphantom{-}2.70\times 10^{-3}$ & $1.57\times 10^{1}$ & $5.20\times 10^{-1}$ & 20 & $3.00\times 10^{-2}$\\
& $\hphantom{-}2.71\times 10^{-3}$ & $1.57\times 10^{1}$ & $1.52\times 10^{-4}$ & 37 & $3.00\times 10^{-2}$\\
\midrule
\multicolumn{6}{c}{Poisson32k}\\
\multirow{3}{*}{PMR\_SPD} 
& $-2.59\times 10^{-2}$ & $1.62\times 10^{1}$ & $3.87\times 10^{1}$ & 5 & $3.09\times 10^{-3}$\\
& $\hphantom{-}1.09\times 10^{-2}$ & $5.92\times 10^{1}$ & $1.37\times 10^{1}$ & 20 & $3.00\times 10^{-2}$\\
& $\hphantom{-}1.03\times 10^{-2}$ & $1.81\times 10^{2}$ & $6.76\times 10^{0}$ & 100 & $3.00\times 10^{-2}$\\
\cmidrule{2-6}
\multirow{3}{*}{PCG} 
& $-2.46\times 10^{-1}$ & $3.97\times 10^{1}$ & $8.70\times 10^{1}$ & 5 & $3.09\times 10^{-3}$\\
& $-7.68\times 10^{-2}$ & $2.94\times 10^{2}$ & $2.89\times 10^{1}$ & 20 & $3.00\times 10^{-2}$\\
& $\hphantom{-}9.92\times 10^{-3}$ & $3.52\times 10^{2}$ & $2.39\times 10^{1}$ & 100 & $3.00\times 10^{-2}$\\
\cmidrule{2-6}
\multirow{3}{*}{LOPMR\_SPD} 
& $-1.98\times 10^{-2}$ & $2.41\times 10^{1}$ & $2.85\times 10^{1}$ & 5 & $3.09\times 10^{-3}$\\
& $\hphantom{-}1.03\times 10^{-2}$ & $2.07\times 10^{2}$ & $6.82\times 10^{0}$ & 20 & $3.00\times 10^{-2}$\\
& $\hphantom{-}1.05\times 10^{-2}$ & $2.70\times 10^{2}$ & $1.17\times 10^{1}$ & 100 & $3.00\times 10^{-2}$\\
\bottomrule
\end{tabular}
\end{table}

\begin{figure}[ht]
\centering
\includegraphics[scale=.8]{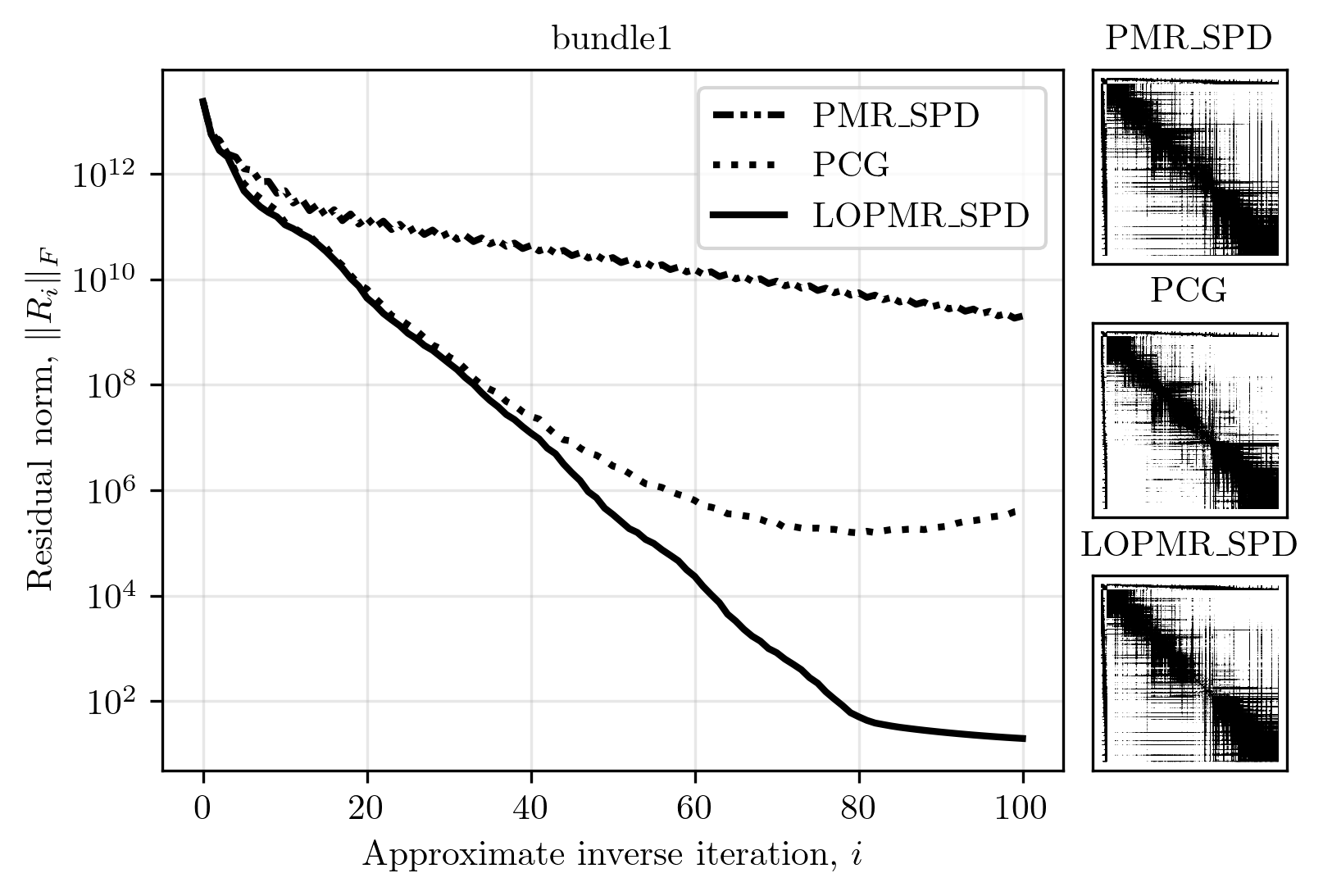}
\includegraphics[scale=.8]{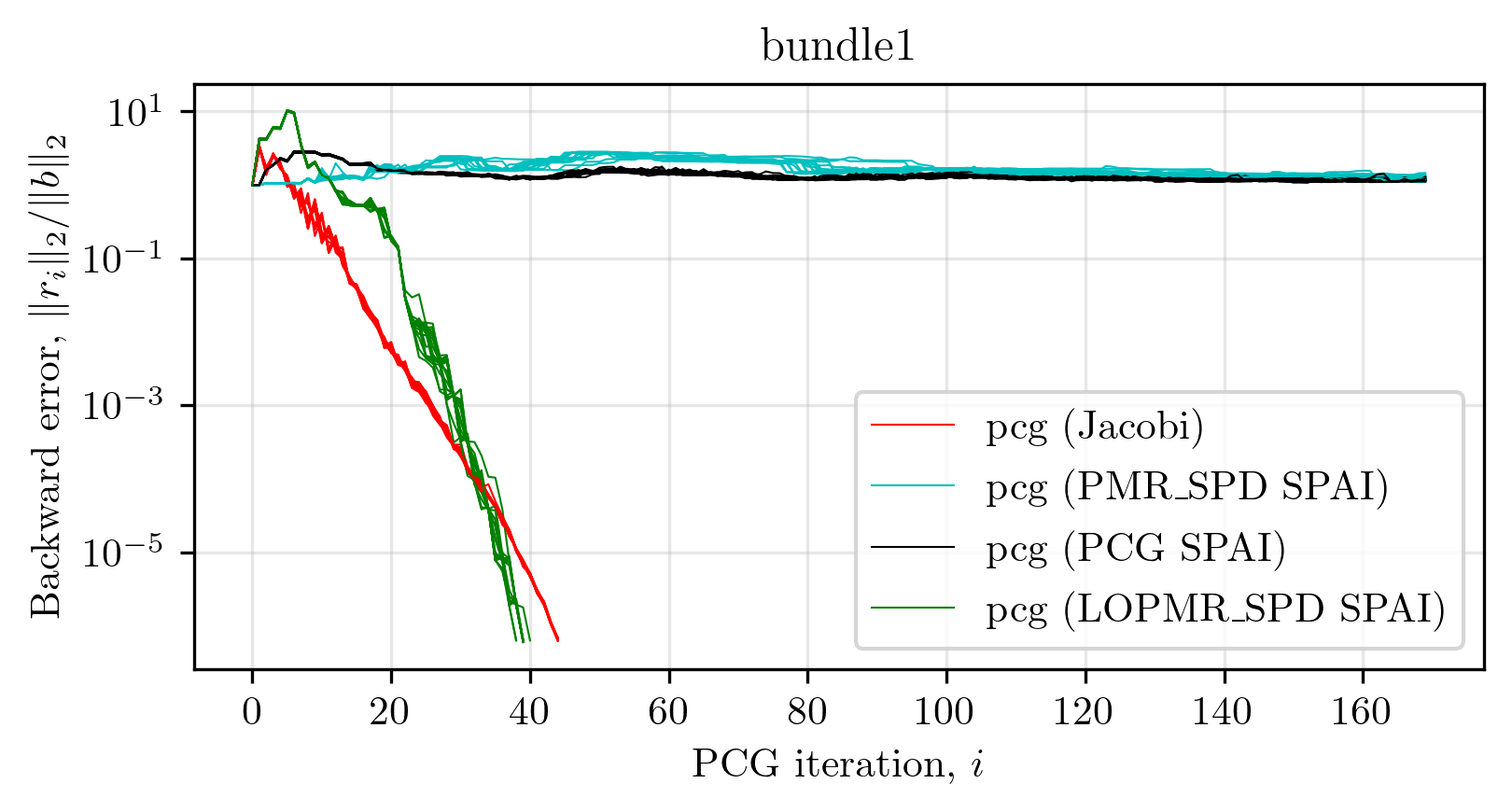}
\caption{Convergence plots and PCG results for SPAIs of the matrix bundle1 obtained with dropping of nonzero values for maximum densities of 3\% (Experiment05-06).}
\label{fig:bundle1}
\end{figure}

The results for the bundle1 matrix, see Figure~\ref{fig:bundle1}, are of particular interest because, as shown in the previous section, the SPAI grows very dense through matrix iterations. 
Consequently, dropping nonzero values becomes the computationally dominant part of matrix iterations, making the additional computations of LOPMR\_SPD and PCG over PMR\_SPD less significant.
As shown in Table~\ref{tab:with-dropping-summary-1}, neither PMR\_SPD nor PCG achieve positive definiteness in 100 iterations or less, whereas LOPMR\_SPD does, even though all the residual Frobenius norms achieved provide no guarantee of non-singularity.
The SPAIs formed by PMR\_SPD and PCG, when used as preconditioners, make vector pcg iterations divergent.
On the other hand, the SPAI formed by LOPMR\_SPD accelerates convergence in comparison to using a Jacobi preconditioner for backward errors below $10^{-4}$.

\begin{figure}[ht]
\centering
\includegraphics[scale=.8]{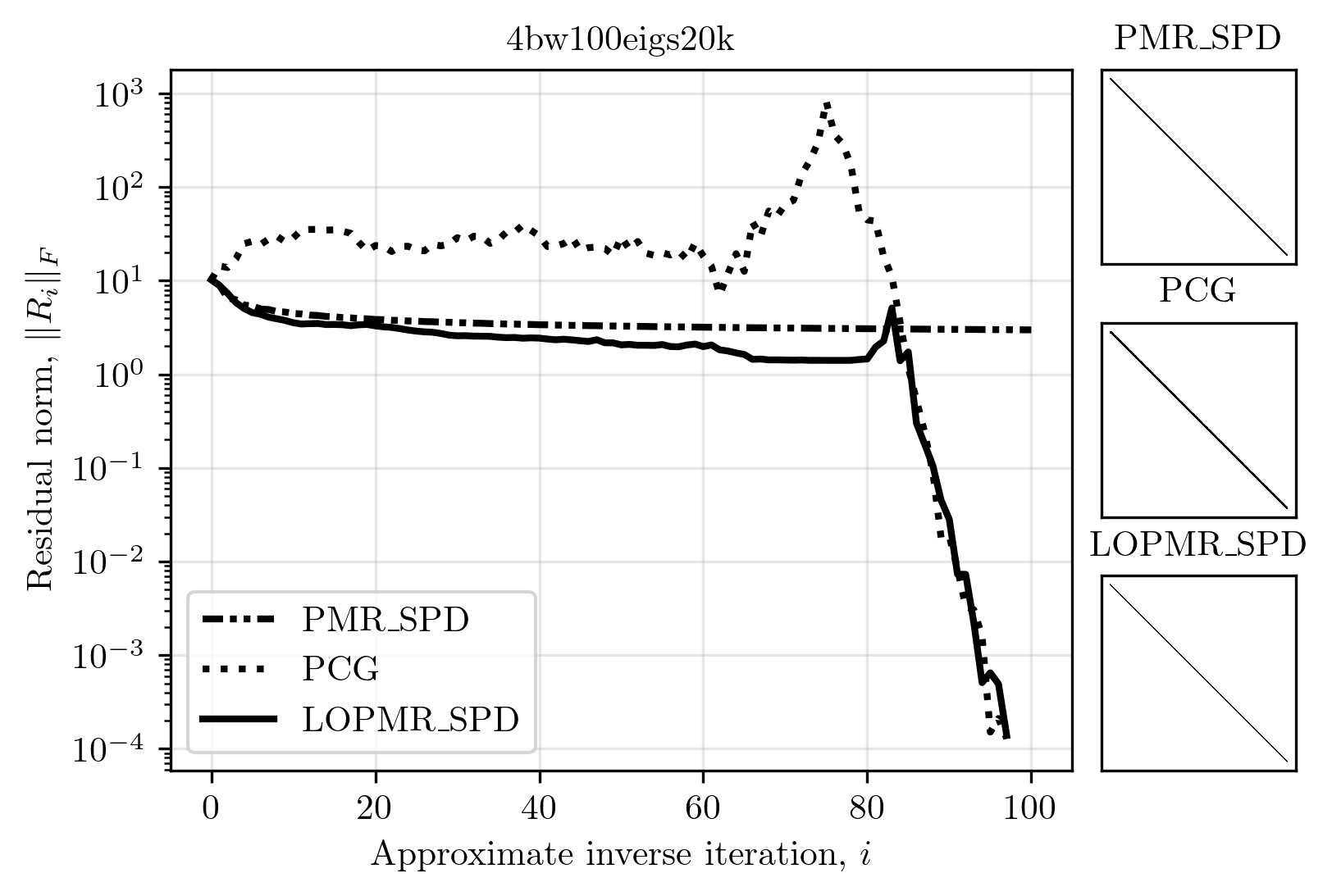}
\includegraphics[scale=.8]{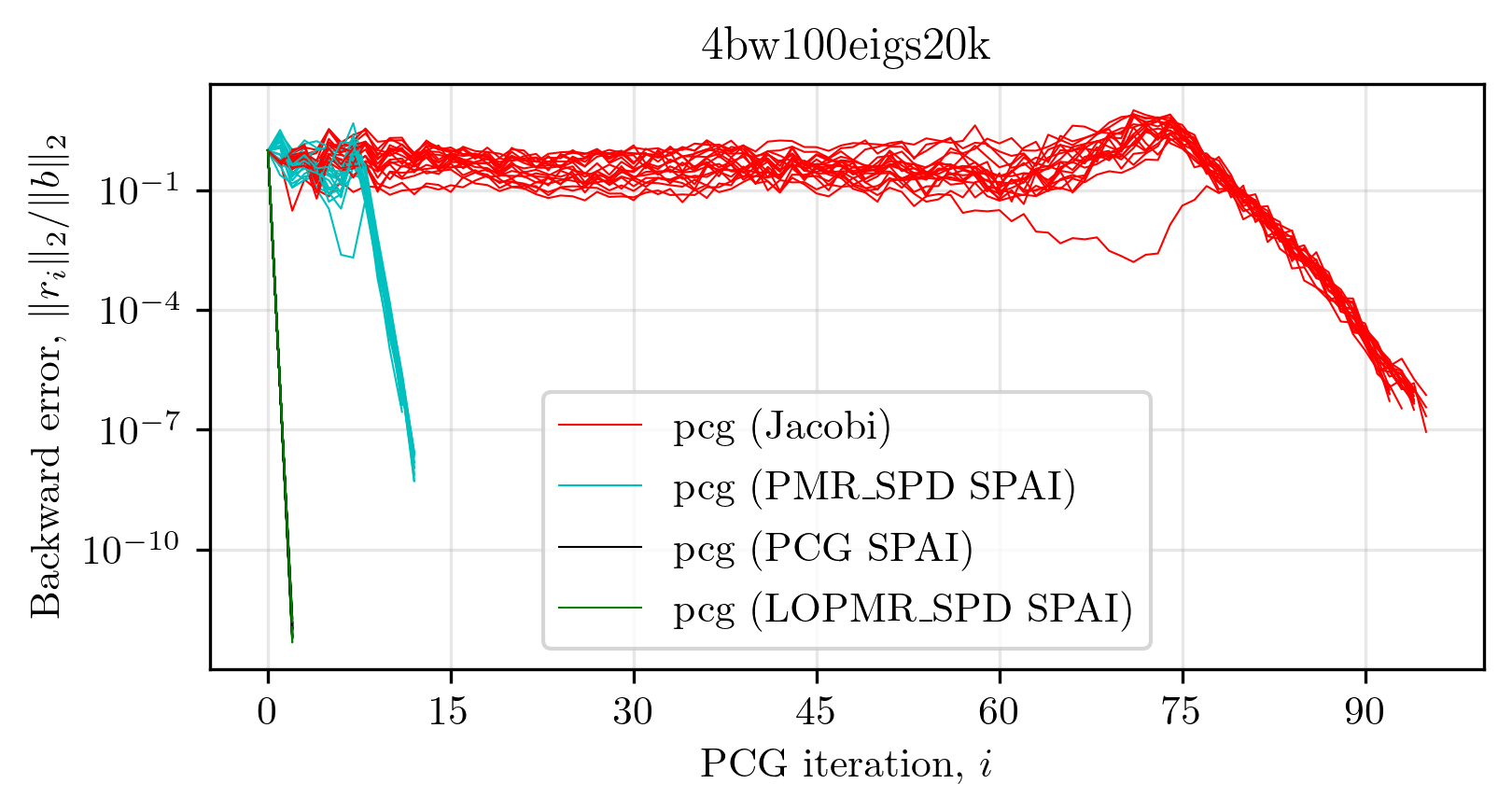}
\caption{Convergence plots and PCG results for SPAIs of the matrix 4bw100eigs20k obtained with dropping of nonzero values for maximum densities of 3\% (Experiment05-06).}
\label{fig:4bw100eigs20k}
\end{figure}

For the 4bw100eigs20k matrix, see Figure~\ref{fig:4bw100eigs20k}, both PCG and LOPMR\_SPD converge to SPAIs of $A$, although the former does so after spurious increases of the matrix residual norm.
After 100 iterations or less, all three methods produce positive definite SPAIs, each of them with less than $0.5\%$ density.
When used as preconditioners of vector pcg iterations, the SPAIs formed by all three methods significantly accelerate convergence in comparison to using a Jacobi preconditioner, but the one of PMR\_SPD less so than those of PCG and LOPMR\_SPD.

\begin{figure}[ht]
\centering
\includegraphics[scale=.8]{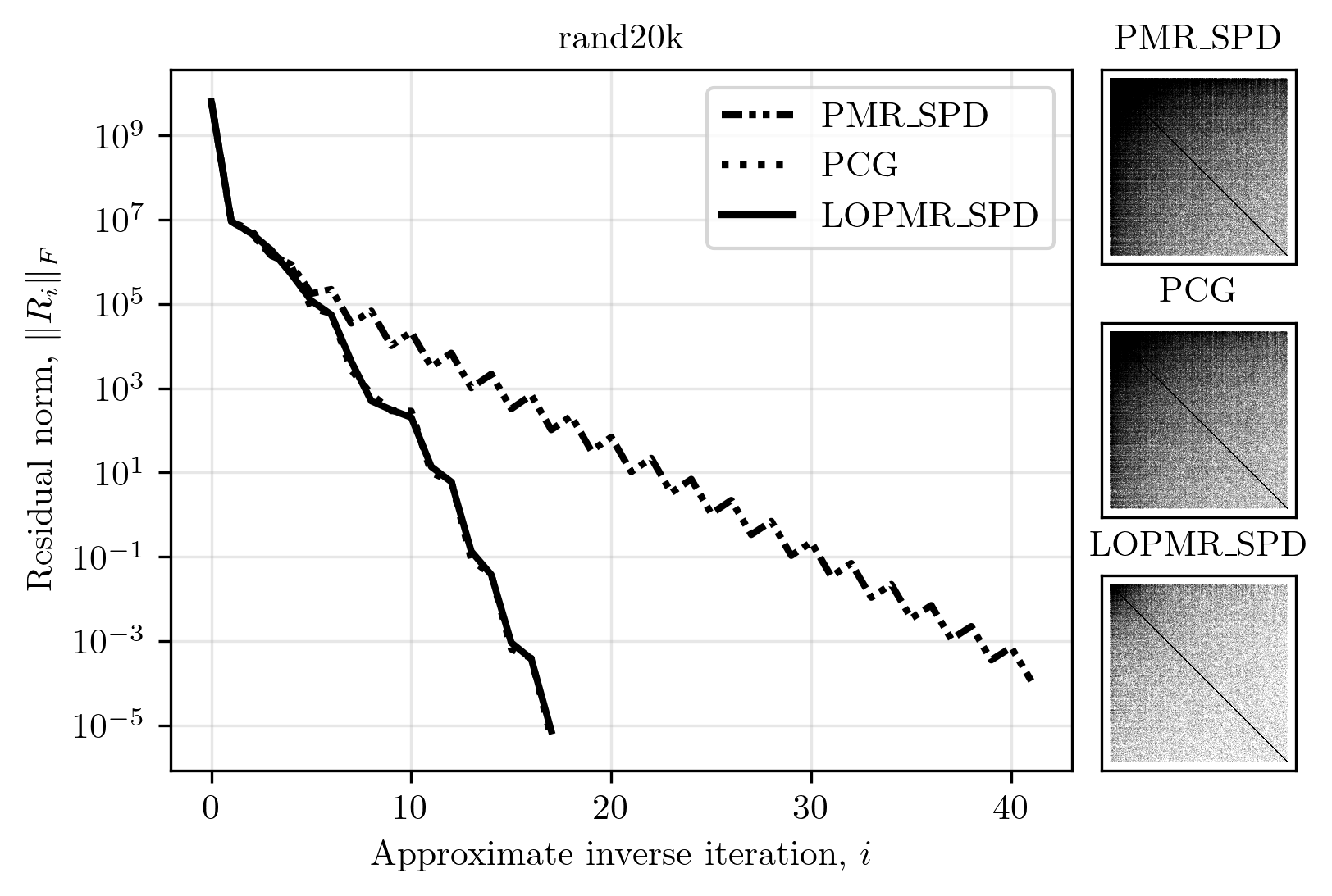}
\includegraphics[scale=.8]{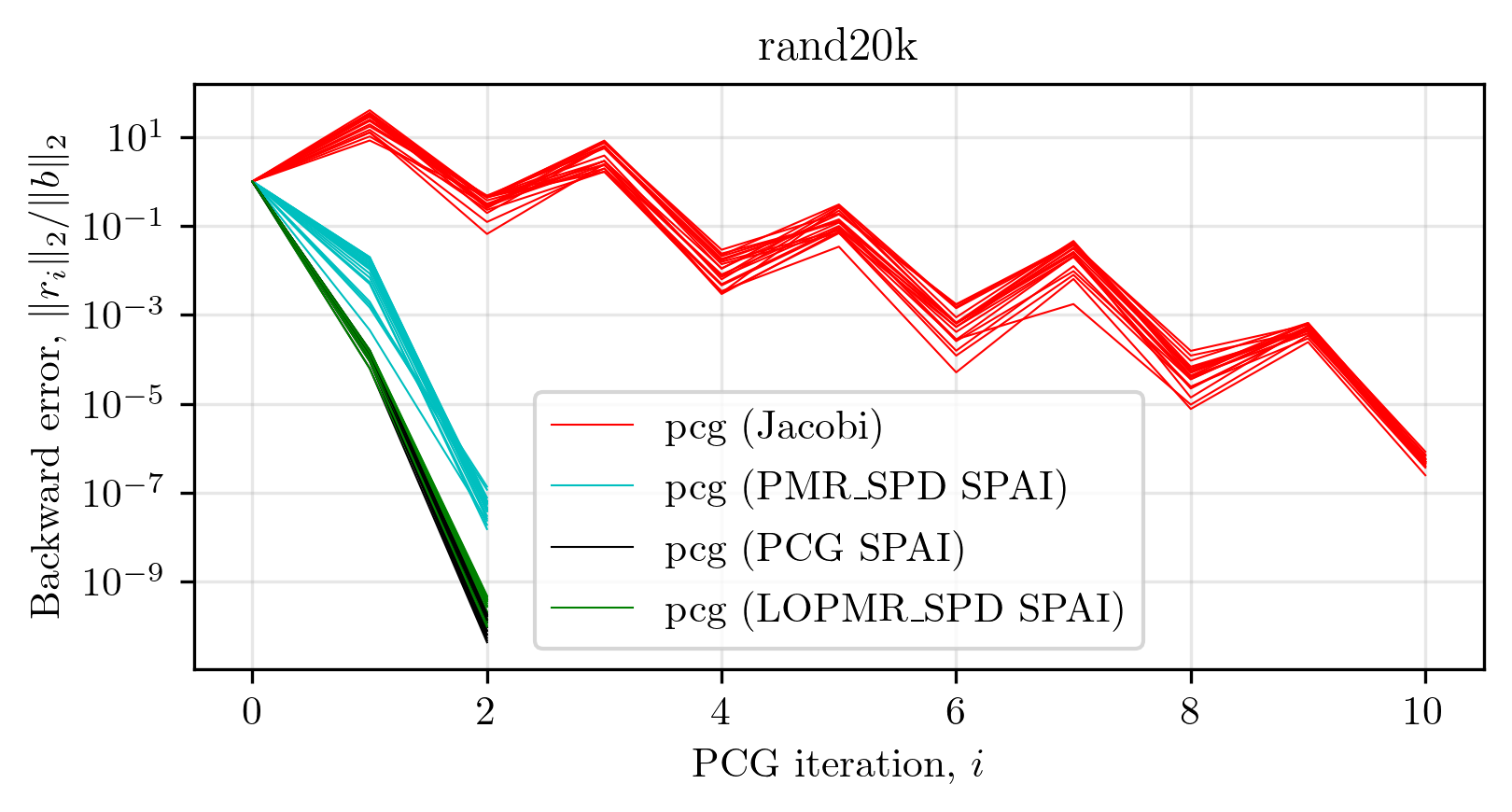}
\caption{Convergence plots and PCG results for SPAIs of the matrix rand20k obtained with dropping of nonzero values for maximum densities of 3\% (Experiment05-06).}
\label{fig:rand20k}
\end{figure}

For the rand20k matrix, see Figure~\ref{fig:rand20k}, all three methods converge to SPAIs of $A$ with backward error below $10^{-6}$.
All the SPAI iterates remain relatively sparse with densities consistently under $0.2\%$.
Interestingly, we note that the combination of matrix iteration with hard thresholding yields a consistent decrease of density through matrix iterations for all three methods.
All three converged SPAIs bring significant convergence acceleration compared to a Jacobi preconditioner when used to precondition vector iteration.

\begin{figure}[ht]
\centering
\includegraphics[scale=.8]{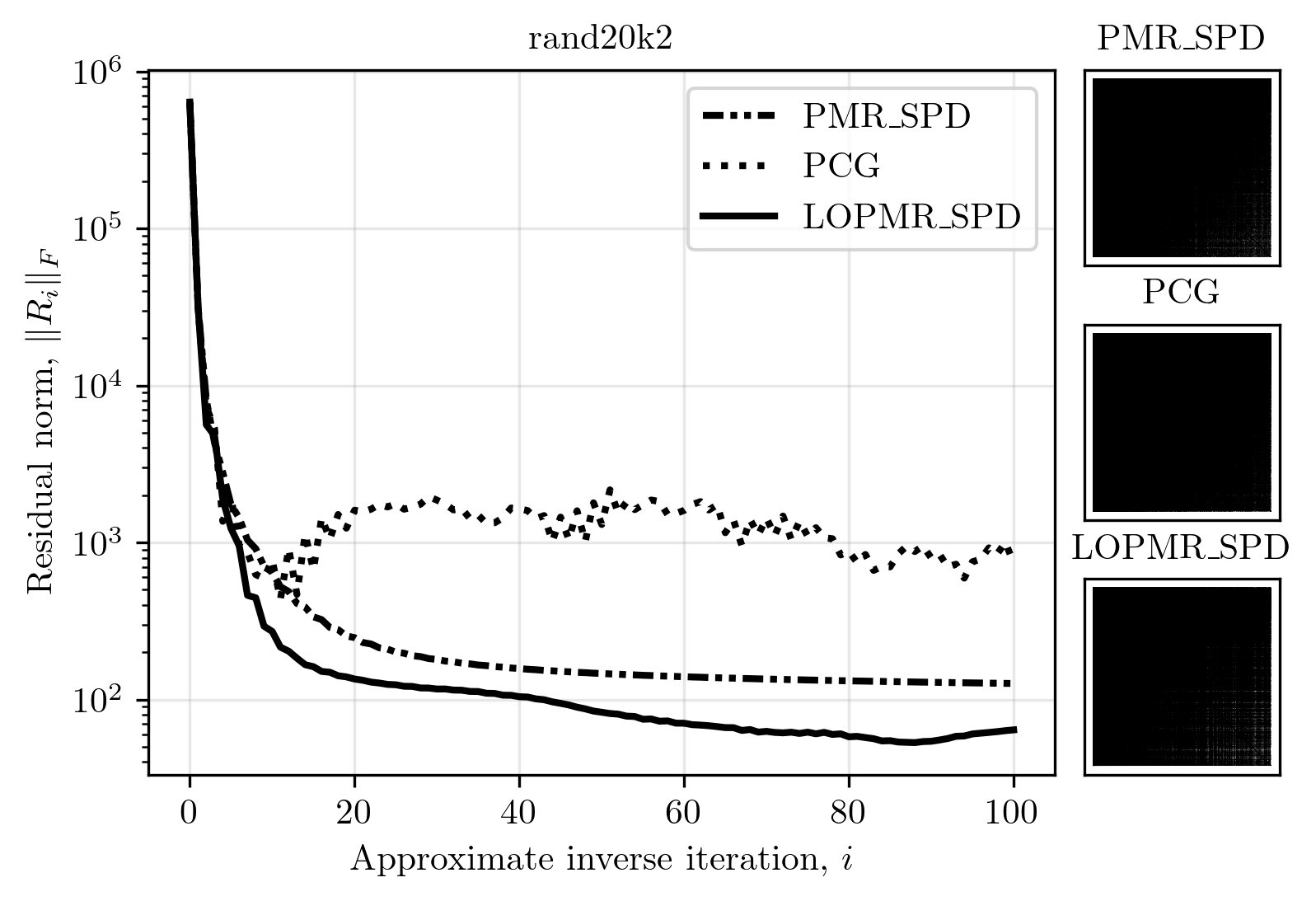}
\includegraphics[scale=.8]{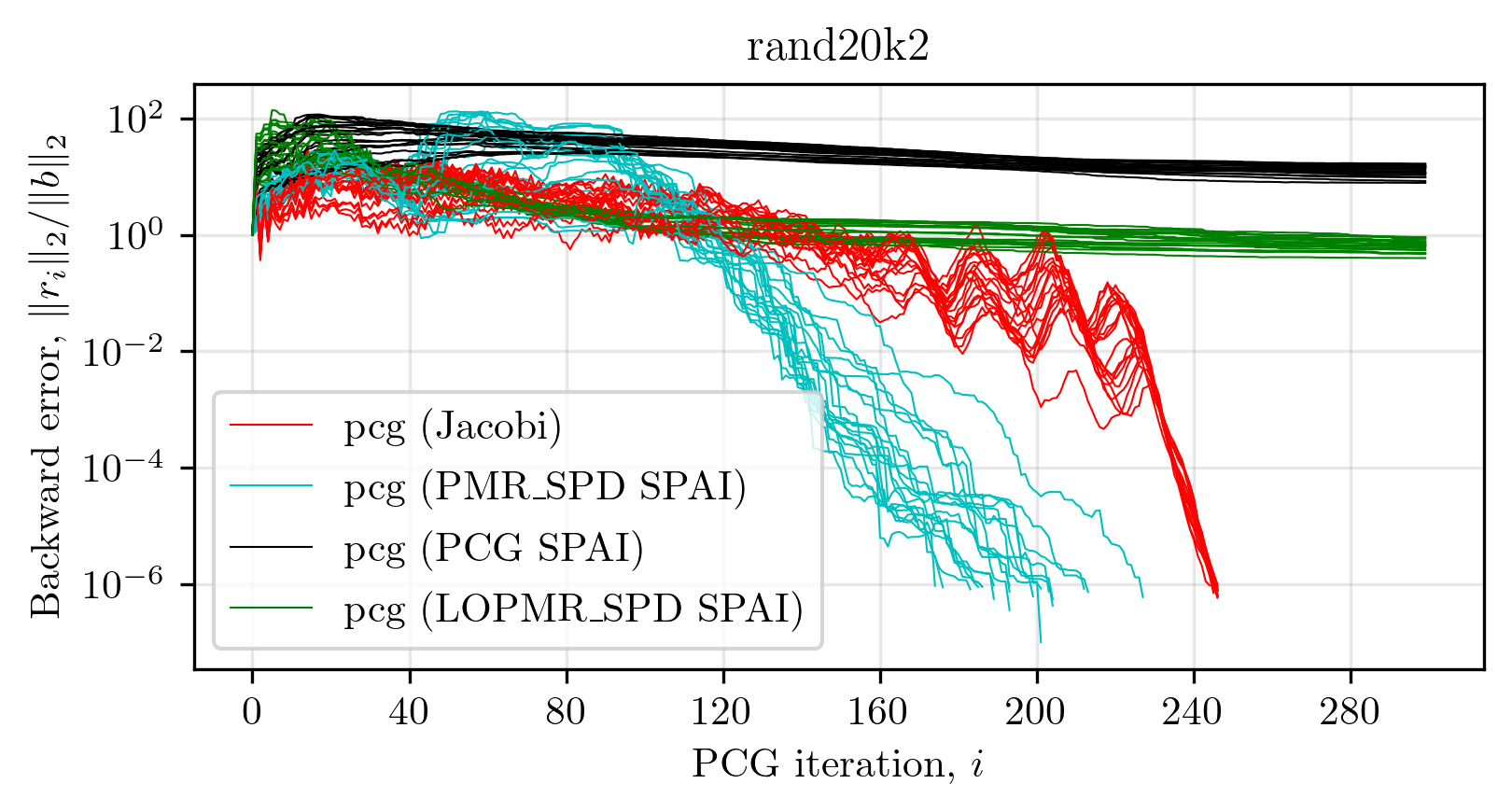}
\caption{Convergence plots and PCG results for SPAIs of the matrix rand20k2 obtained with dropping of nonzero values for maximum densities of 3\% (Experiment05-06).}
\label{fig:rand20k2}
\end{figure}

For the rand20k2 matrix, see Figure~\ref{fig:rand20k2}, despite having the same sparsity pattern as the rand20k matrix, finding a good SPAI proves significantly more difficult as, in 100 matrix iterations, while most iterates 
of PCG and LOPMR\_SPD reach the prescribed density of $3\%$, none converges in residual Frobenius norm.
On the other hand, PMR\_SPD matrix iterates show decreasing densities while not converging either.
Both PMR\_SPD and LOPMR\_SPD matrix iterates do produce positive definite iterates, but not PCG.
Then, when used as preconditioners for vector pcg iteration the surprising result is that the SPAI formed by PMR\_SPD is the only preconditioner accelerating convergence compared to using a Jacobi preconditioner, whereas both PCG and LOPMR\_SPD lead to divergent vector pcg iterations.
%

\begin{figure}[ht]
\centering
\includegraphics[scale=.8]{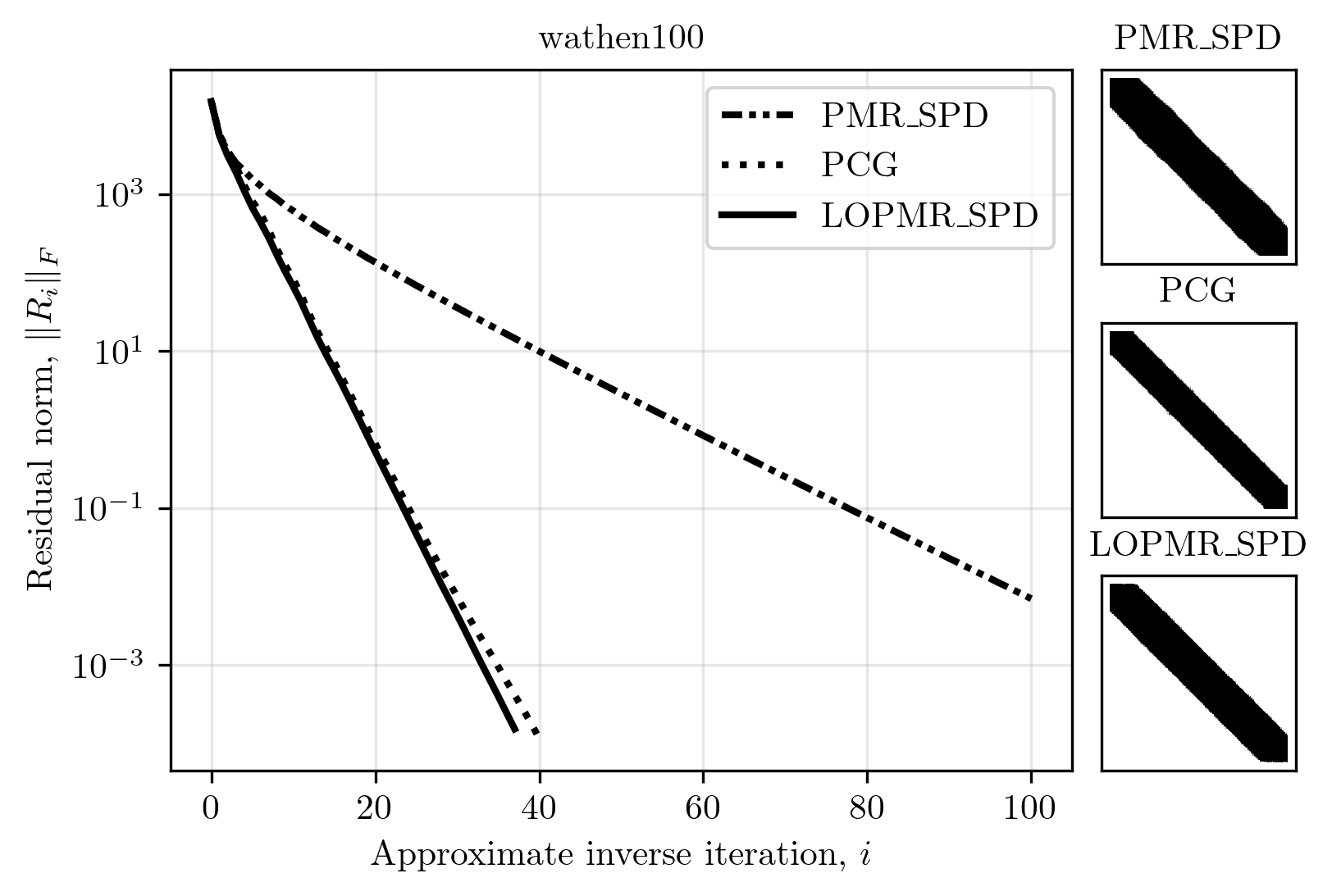}
\includegraphics[scale=.8]{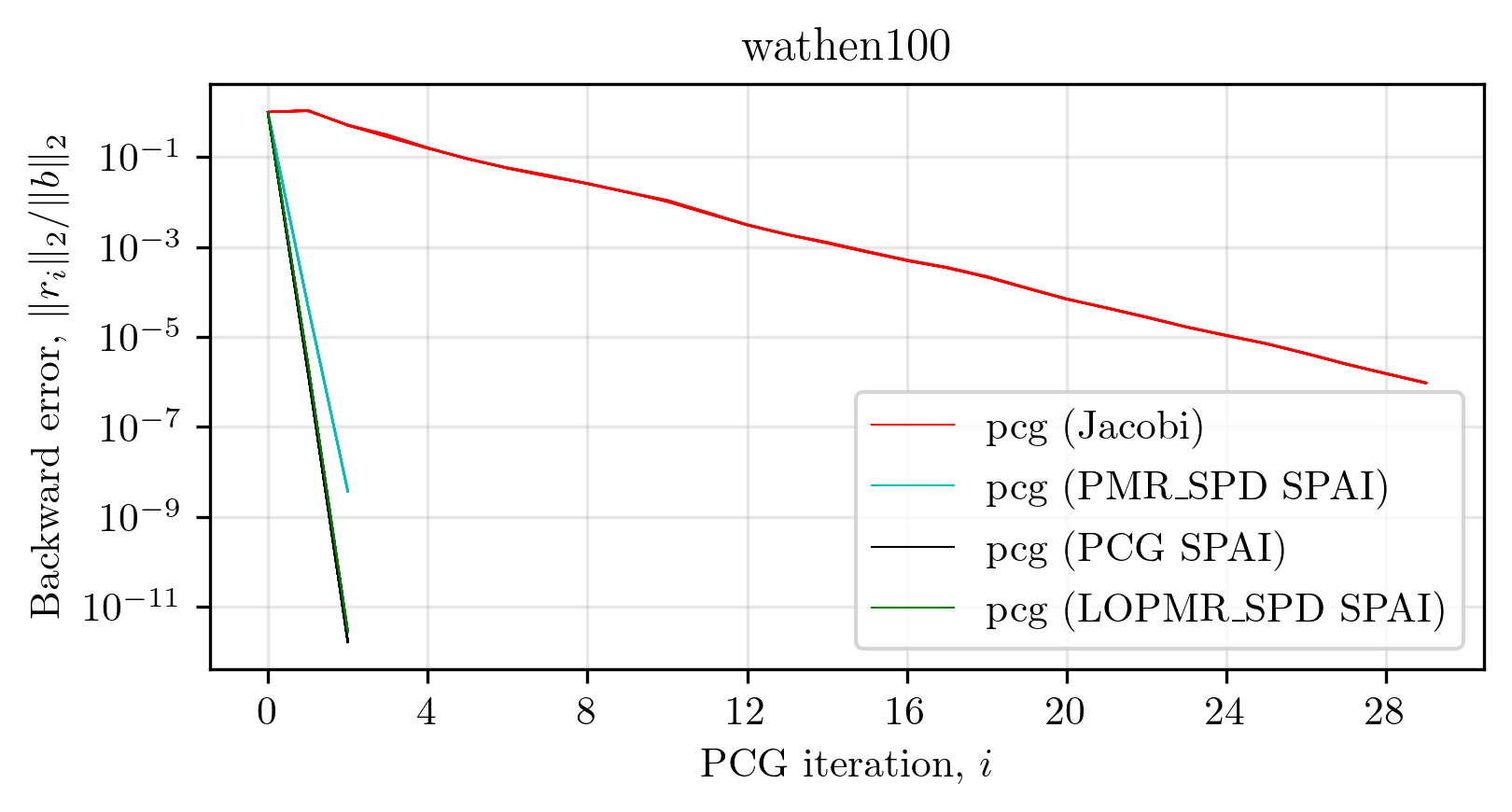}
\caption{Convergence plots and PCG results for SPAIs of the matrix wathen100 obtained with dropping of nonzero values for maximum densities of 3\% (Experiment05-06).}
\label{fig:wathen100}
\end{figure}

\begin{figure}[ht]
\centering
\includegraphics[scale=.8]{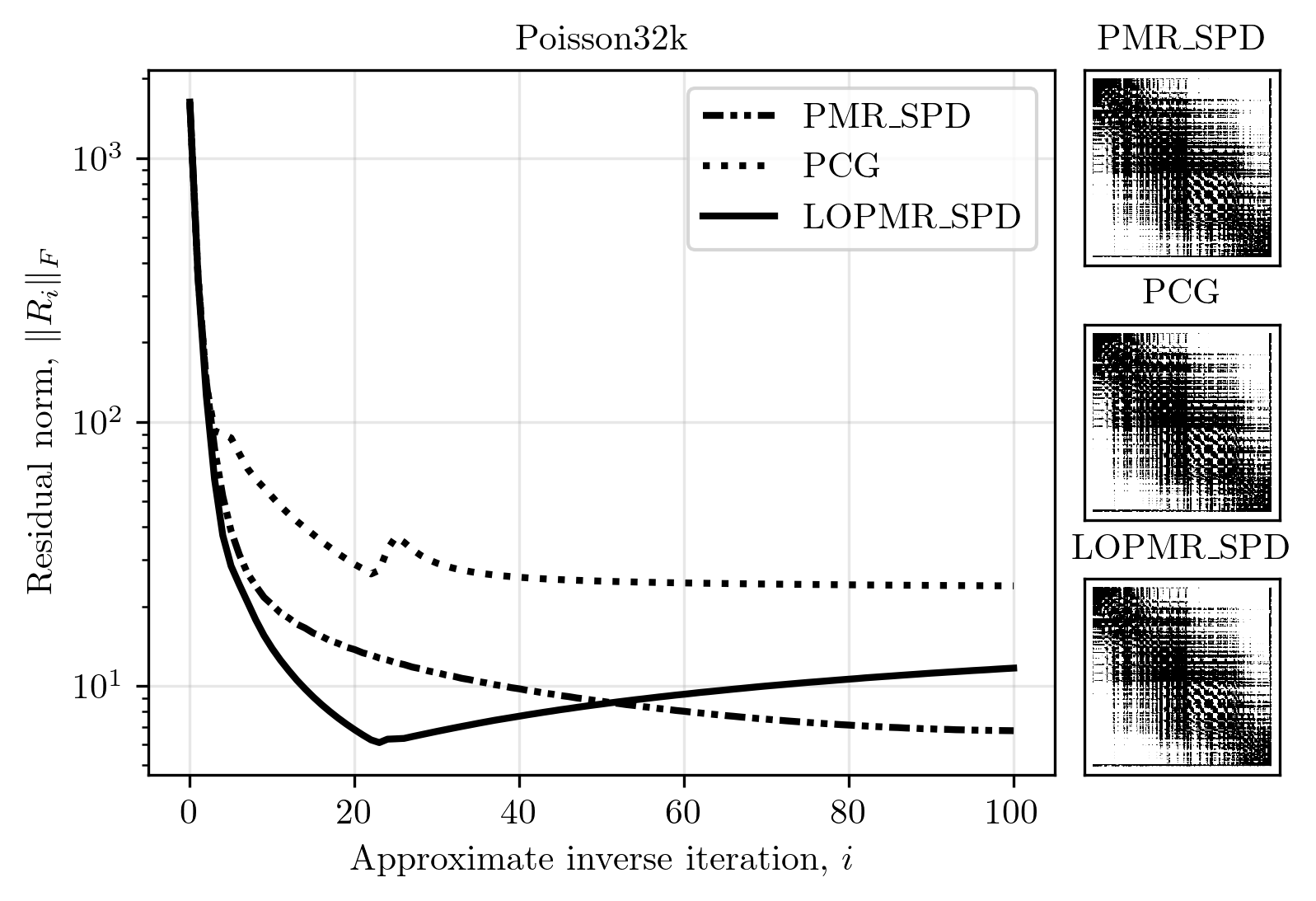}
\includegraphics[scale=.8]{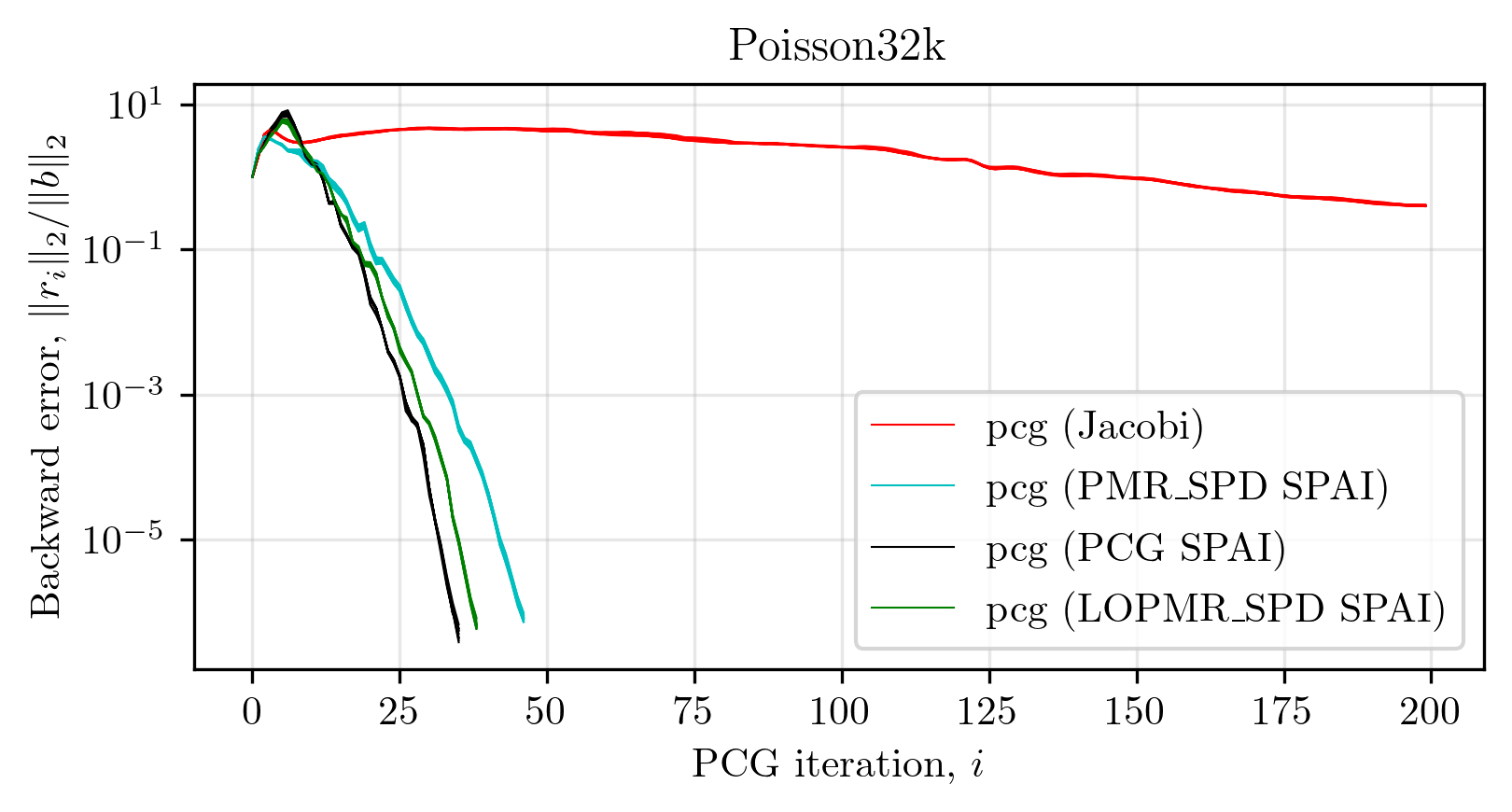}
\caption{Convergence plots and PCG results for SPAIs of the matrix Poisson32k obtained with dropping of nonzero values for maximum densities of 3\% (Experiment05-06).}
\label{fig:Poisson32k}
\end{figure}

For the wathen100 and Poisson32k matrices, see Figures~\ref{fig:wathen100}-\ref{fig:Poisson32k}, although all three methods fail to yield SPAIs with backward error below $10^{-6}$, they do yield preconditioners that significantly accelerate convergence of vector pcg iterations.

In our search for matrices to test the methods presented in this paper, we encountered some matrices for which no good SPD SPAI could be obtained with a maximum density of 3\%, regardless of the number of iterations or method used. 
These matrices are:
\begin{itemize}
\item[-] bodyy6\footnotemark[2]: Moderately conditioned matrix from a structural mechanics problem. Submitted by Alex Pothen from NASA.
\item[-] gyro\_k\footnotemark[2]: Ill-conditioned stiffness matrix from a model order reduction problem.
\item[-] gyro\_m\footnotemark[2]: Moderately conditioned and nearly singular mass matrix from a model order reduction problem.
\item[-] Pres\_Poisson\footnotemark[2]: Moderately conditioned matrix from computational fluid dynamics.
\end{itemize}
\footnotetext[2]{\url{https://sparse.tamu.edu/}}

\section{Conclusion}\label{sec:conclusion}

This work focused on short-recurrence matrix iterations to approximate inverses of sparse SPD matrices.
We presented the SD and MR methods and proved that, for SPD matrices, MR matrix iterations converge linearly to the inverse of $A$ with a rate $(1-\lambda_{\min}(A)^2/\operatorname{tr}(A^2))^{1/2}$, irrespective of the initial iterate. 
LOMR and its preconditioned variant for SPD matrix and preconditioner are introduced as a means to improve the convergence behavior of MR matrix iterations.
LOMR inherits its unconditional convergence from being an enriched variant of MR iterations but, in practice, its convergence behavior is often that of a monotonically decreasing lower envelope to the residual norm of matrix CG iterations, without the occasional spurious oscillations proper to CG matrix iterations.



As evidenced by experiments and theory, the LOMR and CG methods do generally accelerate the convergence behavior in comparison to the SD and MR methods.
The experiments we performed reveal that the LOMR method performs best among all four tested methods, though the CG method exhibits closely matching convergence behavior, albeit with non-monotonic decreases of residual norms and significant oscillations at times.
The success of the LOMR method comes at a cost, namely one SpGEMM, two Frobenius inner products, one Frobenius norm, and one SpGEAM more than CG, the second-best performing method in this work.
Therefore, if CG is capable of achieving a good SPAI, we recommend it over LOMR.
However, as shown through experiments, there are SPD matrices for which only LOMR is capable of achieving a sufficiently good SPAI to serve as a preconditioner.

Some questions remain to be answered.
First, while LOMR and CG prove better than the state-of-the-art MR at building good SPAIs of SPD matrices, their iterates are still not guaranteed to be SPD.
A logical next step would be to develop locally optimal methods analogous to LOMR for the computation of FSAIs for which positive definiteness is consistently enforced for each iterate.
Second, although not documented in this work, we did deploy efforts to try to design effective stopping criteria that guarantee the quality of an SPAI as a preconditioner. 
Unfortunately, those efforts were unsuccessful, and as of now, we see no better approach than to test each iterate of an SPAI as a preconditioner, as in Experiment07.
Third, the methods introduced in this work need to be adapted to parallel computing.
Lastly, a similar study has been initiated for the case of general matrices.

\bibliographystyle{plain} 
\bibliography{references}

\appendix

\section{Dropping heuristic for main iterate}\label{sec:heuristic}
Let the nonzero components of $M$ be denoted by $m_{k\ell}$, and the corresponding residual be $R:=I_n-AM$.
We wish to set multiple nonzero components of $M$ to zero, simultaneously.
In particular, let us introduce the set
\begin{align}\label{eq:400}
\mathcal{D}_M\subset\operatorname{supp}(M):=\{(k,\ell)\text{ such that }m_{k\ell}\neq 0\}
\end{align}
of indices corresponding to the entries we wish to drop, where $\operatorname{supp}(M)$ is the set of indices of nonzero components in $M$.
The perturbed right-approximate inverse iterate is then given by:
\begin{align}\label{eq:405}
\widehat{M}:=M-\sum_{(k,\ell)\in \mathcal{D}_M}m_{k\ell}e_ke_\ell^T
\end{align}
and the corresponding perturbed residual is given by:
\begin{align}\label{eq:410}
\widehat{R}
=R+
\sum_{(k,\ell)\in \mathcal{D}_M}m_{k\ell}Ae_ke_\ell^T.
\end{align}
Then, the Frobenius norm of the perturbed residual becomes
\begin{align}\label{eq:415}
\|\widehat{R}\|_F^2
=&\,\|R\|_F^2
+\sum_{(k,\ell)\in \mathcal{D}_M}\sum_{(r,s)\in \mathcal{D}_M}
m_{k\ell}m_{rs}
\left(Ae_ke_\ell^T,Ae_re_s^T\right)_F
+2\sum_{(k,\ell)\in \mathcal{D}_M} m_{k\ell}\cdot(AR)_{k\ell}
\end{align}
where $(AR)_{k\ell}$ denotes the $(k,\ell)$-entry of $AR$.
Note that the coupling term $(Ae_ke_\ell^T,Ae_re_s^T)_F$ cancels for all $\ell\neq s$.
This means that the effect of nonzero dropping interactions with other droppings only occurs within columns:
\begin{align}\label{eq:420}
\sum_{(k,\ell)\in \mathcal{D}_M}\sum_{(r,s)\in \mathcal{D}_M}
m_{k\ell}m_{rs}
\left(Ae_ke_\ell^T,Ae_re_s^T\right)_F
=
\sum_{t=1}^n
\sum_{\underset{\mathlarger{\ell=t}}{(k,\ell)\in \mathcal{D}_M}}
\sum_{\underset{\mathlarger{s=t}}{(r,s)\in \mathcal{D}_M}}
m_{kt}m_{rt}(Ae_k)^T(Ae_r).
\end{align}
The overall dropping effect on the residual Frobenius norm can then be recast into
\begin{align}\label{eq:425}
\|\widehat{R}\|_F^2-\|R\|_F^2
=&\,
\sum_{(k,\ell)\in \mathcal{D}_M}
m_{k\ell}^2\|Ae_k\|_2^2
+2\sum_{(k,\ell)\in \mathcal{D}_M} m_{k\ell}\cdot(AR)_{k\ell}\\
&\,+\sum_{t=1}^n
\sum_{\underset{\mathlarger{\ell=t}}{(k,\ell)\in \mathcal{D}_M}}
\sum_{\underset{\mathlarger{r\neq k}}{\underset{\mathlarger{s=t}}{(r,s)\in \mathcal{D}_M}}}
m_{kt}m_{rt}(Ae_k)^T(Ae_r).\nonumber
\end{align}
Note that, for the case in which there is only one nonzero component dropped per column, the last term vanishes, in which case we have:
\begin{align}\label{eq:430}
\|\widehat{R}\|_F^2-\|R\|_F^2
=&\,
\sum_{(k,\ell)\in \mathcal{D}_M}
m_{k\ell}^2\|Ae_k\|_2^2
+2\sum_{(k,\ell)\in \mathcal{D}_M} m_{k\ell}\cdot(AR)_{k\ell}.
\end{align}
In practice, we do want to drop more than one nonzero component per column at once, but for practical reasons, we ignore the last term on the right-hand-side of Eq.~\eqref{eq:425}, even when we do drop more than one nonzero element per column.
That is, Eq.~\eqref{eq:430} provides a blueprint to pick nonzero entries of $M^{-1}$ whose dropping leads to the most significant decrease $\|\widehat{R}\|_F^2-\|R\|_F^2$ of residual Frobenius norm.
Although we expect the magnitude of the coupling term ignored in Eq.~\eqref{eq:425} to grow with $|\mathcal{D}_M|$, we find no approach to account for this term in the design of a sufficiently efficient dropping strategy.
Our dropping strategy is as follows:
\begin{enumerate}
\item Symmetrize $M$: $M:=\left(M+M^T\right)/2$.
\item Drop insignificant non-diagonal components: $m_{k\ell}:=0$ for all $|m_{k\ell}|<u$ with $k\neq \ell$, where $u$ denotes the unit round-off. 
\item If a threshold density is achieved, i.e., if $|\operatorname{supp}(M)|>m$ for some fixed $m\ll n^2$, apply the strategy of Eq.~\eqref{eq:435}.
That is,  whenever $|\operatorname{supp}(M)|>m$, we deploy the following strategy:
\end{enumerate}
\begin{align}\label{eq:435}
&\text{For }\;(k,\ell)\in
\mathcal{D}_M=
\arg
\underset{\underset{\mathlarger{|\operatorname{supp}(M)|-|\mathcal{S}|=m}}{\mathcal{S}\subset \operatorname{supp}(M)\setminus\{(k,k),\,k=1,\dots,n\}}}{\min}
\left\{
\sum_{(k,\ell)\in \mathcal{S}}
m_{k\ell}\|Ae_k\|_2^2
+2\sum_{(k,\ell)\in \mathcal{S}} m_{k\ell}\cdot(AR)_{k\ell}
\right\},\\
&\text{ set }\;
m_{k\ell}:=0.\nonumber
\end{align}
Conveniently, the matrix $AR$ is already formed at each iteration in all three non-preconditioned methods presented here.
For the preconditioned methods, extra storage may be needed to store $AR$ along $AZ$.
The dot products $\|Ae_1\|_2^2,\dots,\|Ae_n\|_2^2$ need be computed only once, at the start of the algorithm.

\end{document}